\documentclass[12pt]{article}
\usepackage{epsfig}

\newenvironment{proof}{{\it Proof:\/}}{$\Box$\vskip 0.08in}
\usepackage{amssymb}

\newtheorem{theorem}{Theorem}[section]
\newtheorem{lemma}[theorem]{Lemma}
\newtheorem{corollary}[theorem]{Corollary}

\newtheorem{remark}[theorem]{Remark}

\newtheorem{definition}[theorem]{Definition}
\newtheorem{example}[theorem]{Example}
\newtheorem{property}[theorem]{Property}

\newcommand{\Z}{\mathbb Z}
\newcommand{\Q}{\mathbb Q}

\begin{document}
\begin{center}
\begin{LARGE}
Khovanov homology: torsion and thickness \end{LARGE} \end{center}
\vspace{.2in}

\centerline{\footnotesize Marta M.~Asaeda\footnote{Partially sponsored by 
the NSF grant \#DMS 0202613.}}
\baselineskip=12pt
\centerline{\footnotesize\it Dept. of Mathematics, University of Iowa,} 
\baselineskip=10pt
\centerline{\footnotesize\it 14 MacLean Hall, Iowa City, IA 52242}
\centerline{\footnotesize\it asaeda@math.uiowa.edu}

\vspace*{10pt}
\centerline{\footnotesize  J\'ozef H.~Przytycki } 
\baselineskip=12pt
\centerline{\footnotesize\it Dept. of Mathematics, The George Washington 
University}
\baselineskip=10pt
\centerline{\footnotesize\it 2201 G St. NW, 
Washington, DC 20052}
\centerline{\footnotesize\it przytyck@gwu.edu}

\ \vspace*{.2in}\\
\centerline{e-print July 1, 2004}

{\footnotesize{\begin{minipage}{4.5in} 
\textsc{Abstract.} 
We partially solve the conjecture by A.Shumakovitch that
the Khovanov homology of a prime, non-split link in $S^3$ has 
a non-trivial torsion part.
We give a size restriction  on the Khovanov homology 
of almost alternating links. We relate the Khovanov homology 
of the connected sum of a link diagram with the  Hopf link to 
the Khovanov homology of the diagram via a short exact sequence 
of homology and prove that this sequence splits. 
Finally, we show that our results can be adapted to 
reduced Khovanov homology and that there is a long exact sequence 
connecting reduced Khovanov homology with unreduced homology. 

\end{minipage} }}
\ \\
\ \\

\baselineskip=14pt
\ \\
\begin{Large}{\bf Introduction} \end{Large}

Khovanov homology offers a nontrivial generalization of the Jones 
polynomials of links in $S^3$ (and of the Kauffman bracket skein modules 
of some 3-manifolds). In this paper we use Viro's approach to 
construction of Khovanov homology, and utilize the fact that one works 
with unoriented diagrams (unoriented framed links) in which case there
is a long exact sequence of Khovanov homology. Khovanov homology, over 
the field $\Q$, is a categorification of the Jones polynomial (i.e. we 
represent the Jones polynomial as the 
generating function of Euler characteristics). However, for integral 
coefficients Khovanov homology almost always has torsion. The first 
part of the paper  is devoted to the construction of torsion in 
Khovanov homology. In the second part of the paper we analyze the 
thickness of Khovanov homology and reduced Khovanov homology.

The paper is organized as follows.\ 
In the first section we recall the definition of Khovanov homology and
its basic properties. 

In the second section we prove that adequate link diagrams with 
an odd cycle property have $\Z_2$-torsion in Khovanov homology.

In the third section we discuss torsion in the Khovanov homology of an 
 adequate link diagram with an even cycle property.

In the fourth section we prove Shumakovitch's theorem that 
prime, non-split alternating links different from the trivial knot and 
the Hopf link have $\Z_2$-torsion in Khovanov homology.
We generalize this result to a class of adequate links.

In the fifth section we generalize result of E.S.Lee about 
the Khovanov homology of alternating links (they are 
$H$-thin\footnote{We also 
found a simple proof of Lee's result\cite{Lee-1,Lee-2} 
that for alternating links Khovanov 
homology yields the classical signature, see Remark 1.6.}). 
We do not assume rational coefficients in this 
generalization and we allow alternating adequate links on a surface.
We use Viro's exact sequence of Khovanov homology to extend Lee's 
results to  almost alternating diagrams and $H$-$k$-thick links.

In the sixth section we compute the Khovanov homology for a connected sum of 
$n$ copies of Hopf links and construct a short exact sequence 
of Khovanov homology 
involving a link and its connected sum with the Hopf link. By showing that 
this sequence splits, we answer the question asked by Shumakovitch. 

In the seventh section we notice that the results of sections 5 and 6 can be 
adapted to reduced Khovanov homology. Finally, we show that 
there is a long exact sequence
connecting reduced Khovanov homology with unreduced homology.
 
\section{Basic properties of Khovanov homology}\label{1}

The first spectacular application of the Jones polynomial 
(via Kauffman bracket relation) was the solution of Tait conjectures on 
alternating diagrams and their generalizationsm to adequate diagrams.
Our method of analysing torsion in Khovanov homology has its root in 
work related to solutions of Tait conjectures \cite{Ka,Mu,Thi}. 

Recall that the Kauffman bracket polynomial $<D>$ of a link diagram $D$
is defined by the skein relations
$<\parbox{0.6cm}{\psfig{figure=L+nmaly.eps,height=0.6cm}}> =
A <\parbox{0.6cm}{\psfig{figure=L0nmaly.eps,height=0.6cm}}> + A^{-1}
<\parbox{0.6cm}{\psfig{figure=Linftynmaly.eps,height=0.6cm}}>$
and $<D\sqcup \bigcirc> =(-A^2-A^{-2})<D>$ and the normalization
$<\bigcirc> = 1$. The categorification of this invariant (named
by Khovanov {\it reduced homology}) is discussed in Section 7.
For the (unreduced) Khovanov homology we use the version of 
the Kauffman bracket polynomial
normalized to be $1$ for the empty link (we use the notation $[D]$ 
in this case).

\begin{definition}[Kauffman States]\label{1.1} \ \\
Let $D$ be a diagram\footnote{We think of the 3-ball $B^3$ as $D^2\times I$ 
and the diagram is drawn on the disc $D^2$. In \cite{APS}, 
we have proved that 
the theory of Khovanov homology can be extended to links in an oriented 
3-manifold $M$ that is the bundle over a surface 
$F$ ($M= F \tilde \times I$). 
If $F$ is orientable then
$M= F\times I$. If $F$ is unorientable then $M$ is a twisted 
$I$ bundle over $F$ (denoted by $F \hat \times I$). Several results 
of the paper are valid for the Khovanov homology of links 
in $M= F \tilde \times I$.} of an unoriented, framed link in a 3-ball 
$B^3$.  A Kauffman state $s$
of $D$ is a function from the set of crossings of $D$ to
the set $\{+1,-1\}$. Equivalently,
we assign to each crossing of $D$ a marker
according to the following convention:\\
\ \\
\centerline{\psfig{figure=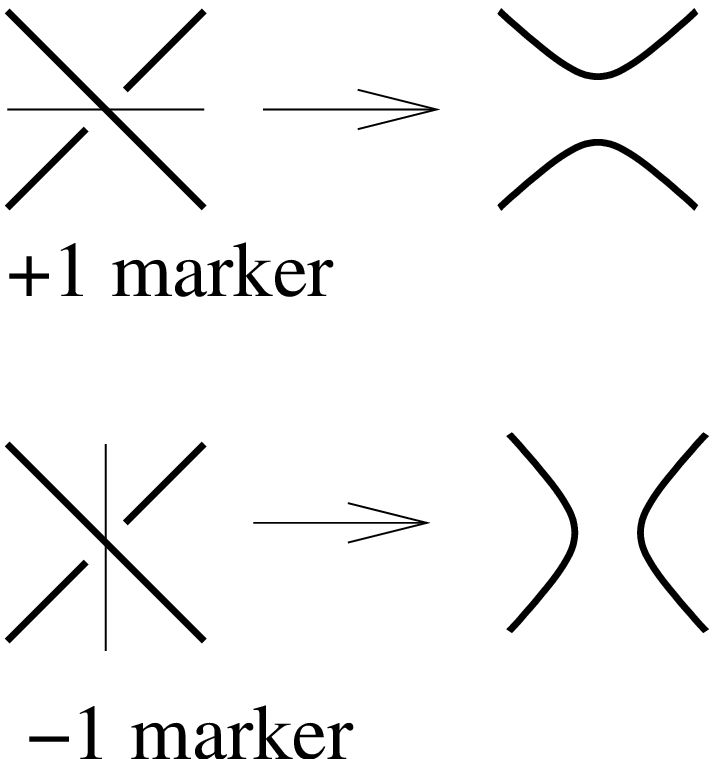,height=5.5cm}}
\begin{center}
Fig. 1.1; markers and associated smoothings
\end{center}

By $D_s$ we denote the system of circles in the diagram 
obtained by smoothing all
crossings of $D$ according to the markers of the state $s$, Fig. 1.1.\\
By $|s|$ we denote the number of components of $D_s$.
\end{definition}

Using this notation we have the Kauffman bracket polynomial given by
the state sum formula:
$[D] = (-A^2-A^{-2})<D> = \sum_s A^{\sigma(s)}(-A^2-A^{-2})^{|s|}$,
where $\sigma(s)$ is the number
of positive markers minus the number of negative markers in the state $s$.


To define Khovanov homology it is convenient (as noticed by Viro) to 
consider enhanced Kauffman states.
\begin{definition}\label{1.2}
An enhanced Kauffman state $S$ of
an unoriented framed link diagram
$D$ is a Kauffman state $s$ with an additional assignment of
$+$ or $-$ sign to each circle of $D_s.$
\end{definition}
Using enhanced states we express the Kauffman bracket polynomial as
a (state) sum of monomials which is important in the definition
of Khovanov homology we use. We have
$[D] = (-A^2-A^{-2})<D> = \sum_S (-1)^{\tau(S)}A^{\sigma(s)+2\tau(S)}$,
where $\tau(S)$ is the number of positive circles minus 
the number of negative circles in the enhanced state $S$.

\begin{definition}[Khovanov chain complex]\label{1.3}\ \
\begin{enumerate}
\item[(i]) Let ${\cal S}(D)$ denote the set of enhanced Kauffman states 
of a diagram $D$, and let ${\cal S}_{i,j}(D)$ denote the set of enhanced 
Kauffman states $S$ such that $\sigma(S) = i$ and $\sigma(S) +2\tau(S) = j$, 
The group ${\cal C}(D)$ 
(resp. ${\cal C}_{i,j}(D)$) is defined to be the free abelian group 
spanned by ${\cal S}(D)$ (resp. ${\cal S}_{i,j}(D)$). ${\cal C}(D) = 
\bigoplus_{i,j\in \Z} {\cal C}_{i,j}(D)$ is a free abelian 
 group with (bi)-gradation.  
\item[(ii]) For a link diagram $D$ with ordered crossings, we define the 
chain complex $({\cal C}(D),d)$ where $d=\{d_{i,j}\}$ and 
the differential $d_{i,j}: {\cal C}_{i,j}(D) \to 
{\cal C}_{i-2,j}(D)$ satisfies $d(S) = \sum_{S'} (-1)^{t(S:S')}[S:S'] S'$ 
with $S\in {\cal S}_{i,j}(D)$, $S'\in {\cal S}_{i-2,j}(D)$, and 
 $[S:S']$  
equal to $0$ or $1$. $[S:S']=1$ if and only if markers of $S$ and $S'$ 
differ exactly at one crossing, call it $c$, and all the circles 
of $D_S$ and $D_{S'}$ 
not touching $c$ have the same sign\footnote{From our conditions 
it follows that at the crossing $c$ the marker of $S$ is positive, 
 the marker of $S'$ is negative, and 
that $\tau(S') = \tau(S)+1$.}. Furthermore, $t(S:S')$ is the number of 
negative markers assigned to crossings in $S$ bigger than $c$ in 
the chosen ordering.
\item[(iii)] The Khovanov homology of the diagram $D$ is defined to be 
the homology of the chain complex $({\cal C}(D),d)$; 
$H_{i,j}(D) = ker(d_{i,j})/d_{i+2,j}({\cal C}_{i+2,j}(D))$. The 
Khovanov cohomology of the diagram $D$ are defined to be the cohomology 
of the chain complex $({\cal C}(D),d)$.
\end{enumerate}
\end{definition}

Below we list a few elementary properties of Khovanov homology
following from properties of Kauffman states used in the proof
of Tait conjectures \cite{Ka,Mu,Thi}.

The positive state $s_+=s_+(D)$ (respectively the negative
state $s_-=s_-(D)$) is the state with all positive markers
(resp. negative markers). The alternating diagrams without nugatory 
crossings (i.e. crossings in a diagram of the form
\parbox{2.2cm}{\psfig{figure=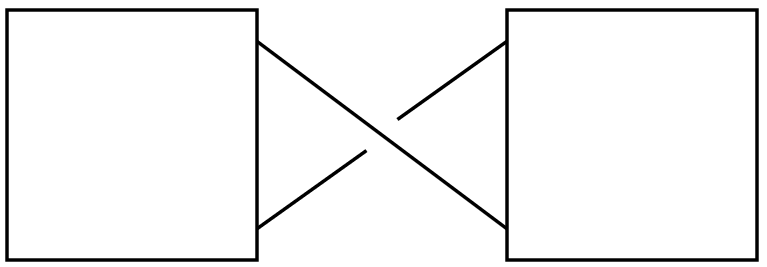,height=0.7cm}})
are generalized 
to adequate diagrams using properties of states $s_+$ and $s_-$. 
Namely, the diagram $D$ is $+$-adequate (resp. $-$-adequate) 
 if the state of positive (resp. negative) markers,
$s_+$ (resp. $s_-$), cuts the diagram to
the collection of circles, so that every crossing is
connecting different circles. $D$ is an adequate diagram if it is 
$+$- and $-$-adequate \cite{L-T}.

\begin{property}\label{1.4}\ \\
If $D$ is a diagram of $n$ crossings and its positive state $s_+$ has
$|s_+|$ circles then the highest term (in both grading indexes) 
 of Khovanov chain complex is
${\cal C}_{n,n+2|s_+|}(D)$; we have ${\cal C}_{n,n+2|s_+|}(D)= \Z$.
 Furthermore, if $D$ is a $+$-adequate diagram, 
 then the whole group
${\cal C}_{*,n+2|s_+|}(D)=\Z$ and $H_{n,n+2|s_+|}(D)= \Z$. 
Similarly the lowest term in the Khovanov
chain complex is ${\cal C}_{-n,-n-2|s_-|}(D)$. 
Furthermore, if $D$ is a $-$-adequate diagram,
 then the whole group
${\cal C}_{*,-n-2|s_-|}(D)=\Z$ and $H_{-n,-n-2|s_-|}(D)= \Z$.
Assume that $D$ is
a non-split diagram then $|s_+| +|s_-| \leq n+2$ and the equality
holds if and only if $D$ is an alternating diagram or
a connected sum of such diagrams (Wu's dual state lemma \cite{Wu}).
\end{property}
\begin{property} \label{1.5} \
Let $\sigma(L)$ be the classical (Trotter-Murasugi) 
signature\footnote{One should not mix the signature $\sigma(L)$ with
$\sigma(s)$ which is the signed sum of markers of the state $s$ of 
a link diagram.} of 
an oriented link $L$ and $\hat \sigma (L) = \sigma({L}) + lk({L})$, 
where $lk(L)$ is the global linking number of $L$, its Murasugi's 
version which does not depend on an orientation of $L$. Then
\begin{enumerate}
\item[(i)] [Traczyk's local property]\ 
If $D^v_0$ is a link diagram obtained from
an oriented alternating link diagram $D$ by smoothing its crossing 
$v$ and $D^v_0$ has the same number of (graph) components as $D$, then
$\sigma(D) = \sigma(D^v_0)- sgn(v)$. 
One defines the sign of a crossing $v$ as\
  $sgn(v) = \pm 1$ according to
the convention $sgn(\parbox{0.9cm}{\psfig{figure=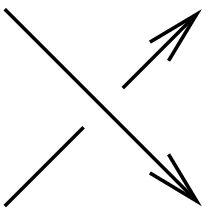,height=0.8cm}}
)= 1$ and $sgn(\parbox{0.9cm}{\psfig{figure=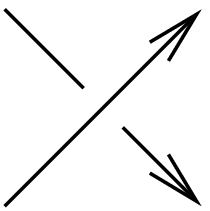,height=0.8cm}}
)=-1$.
\item[(ii)] [Traczyk Theorem \cite{Tra,Pr}] \ \\
The signature, $\sigma (D)$, of the non-split alternating oriented 
link diagram $D$ is equal to 
$n^- - |s_-| +1 = -n^+ + |s_+|-1 = -\frac{1}{2}(n^+ - n^- 
-(|s_+|- |s-|))=$ $n^- - n^+ + d^+ - d^-$, where $n^+(D)$ (resp. $n^-(D)$) 
is the number of positive (resp. negative) crossings
of $D$ and $d^+$ (resp. $d^-$) is the 
number of positive (resp.  negative) edges in a spanning forest of 
the Seifert graph\footnote{The Seifert graph, $GS(D)$, of 
an oriented link diagram $D$ is 
a signed graph whose vertices are in bijection with Seifert circles of $D$ 
and edges are in a natural bijection with crossings of $D$. For an alternating 
diagram the 2-connected components (blocks) of $GS(D)$ have edges of the 
same sign which makes $d^+$ and $d^-$ well defined.}  
of $D$. 
\item[(iii)]
[Murasugi's Theorem \cite{Mu-1,Mu-2}]\ \\
Let $D$ be a non-split alternating oriented diagram without
nugatory crossings or a
connected sum of such diagrams. 
Let $V_D(t)$ be its Jones polynomial\footnote{Recall that if 
$\vec{D}$ is an oriented diagram (any orientation
put on the unoriented diagram $D$),  and $w(\vec{D})$ is
its writhe or Tait number, $w(\vec{D}) = n^+-n^-$, 
then $V_{\vec{D}}(t) = A^{-3w(\vec{D})}<D>$ for $t=A^{-4}$.
P.G.Tait (1831-1901) was the first to consider the number $w(\vec{D})$
 and is often called
the Tait number of the diagram $\vec{D}$ and denoted by
$Tait(\vec{D})$.},
 then the maximal degree $\max \ V_L(t) = n^+(L) - \frac{\sigma(L)}{2}$
and the minimal degree $\min \ V_L(t) = -n^-(L) - \frac{\sigma(L)}{2}$.
\item[(iv)] [Murasugi's Theorem for unoriented link diagrams].
Let $D$ be a non-split alternating unoriented diagram without
nugatory crossings or a connected sum of such diagrams. \\
Then the maximal degree 
$\max \ <D> = \max \ [D] -2 = n+2|s_+|-2=2n+ sw(D)+2\hat\sigma(D)$\\
and the minimal degree 
$\min \ <D> = \min \ [D] +2 = -n - 2|s_-|+2= -2n+sw(D)+ 2\hat\sigma(D)$. 
The self-twist number of a diagram $sw(D)=
\sum_v sgn(v)$, where the sum is
taken over all self-crossings of $D$. A self-crossing involves arcs
from the same component of a link. $sw(D)$ does not depend on orientation 
of $D$.
\end{enumerate}
\end{property}
\begin{remark}\label{1.6} \ \\ 
In Section 5 we reprove the result of Lee \cite{Lee-1} 
that the Khovanov homology 
of non-split alternating links is supported by two adjacent 
diagonals of slope 2, that is $H_{i,j}(D)$ can be nontrivial only for 
two values of $j-2i$ which differ by $4$ (Corollary 5.5). 
One can combine Murasugi-Traczyk result with Viro's 
long exact sequence of Khovanov  
homology  and Theorem 7.3 to recover Lee's result (\cite{Lee-2}) 
that for alternating links
Khovanov homology has the same information as the Jones polynomial and
the classical signature\footnote{The beautiful
paper by  Jacob Rasmussen \cite{Ras} generalizes Lee's results  
and fulfil our dream
(with Pawe{\l} Traczyk) of constructing a ``supersignature" from Jones type
construction \cite{Pr}.} (see Chapter 10 of \cite{Pr}). 
From properties 1.5 and 1.6 it follows that 
for non-split alternating diagram without nugatory crossings
$H_{n,2n+sw+2\hat\sigma +2}(D)=H_{-n,-2n+sw+2\hat\sigma -2}(D)= \Z$. 
Thus diagonals which support nontrivial $H_{i,j}(D)$ satisfy $j-2i= 
sw(D) +2\hat\sigma(D) \pm 2$. If we consider Khovanov cohomology 
$H^{i',j'}(D)$, as considered in \cite{Kh-1,BN-1}, 
then $H^{i',j'}(D)=H_{i,j}(D)$ 
for $i'=\frac{w(D)-i}{2}$, $j'=\frac{3w(D)-j}{2}$ and thus $j'-2i' = 
\frac{-1}{2}(j-2i-w(D))=\sigma(D)\mp 1$ as in Lee's Theorem.
\end{remark}
\begin{remark}\label{1.7}\ \\
The definition of Khovanov homology extends to links in $I$-bundles over 
surfaces $F$ ($F\neq RP^2$) \cite{APS}. 
In the definition we must differentiate between 
trivial curves, curves bounding a M\"obius band, and other non-trivial curves.
Namely, we define $\tau (S)$ as the sum of signs of circles of $D_S$ taken 
over all trivial circles of $D_S$. Furthermore, 
to have $[S:S']=1$, we assume additionally that the sum of signs 
of circles of $D_S$ taken over all nontrivial circles which do not bound 
a M\"obius band is the same for $S$ and $S'$.
\end{remark}

\section{Diagrams with odd cycle property}\label{2}
In the next few sections we use the concept of a graph, $G_s(D)$,
 associated to a link diagram $D$ and its state $s$. The graphs 
corresponding to states $s_+$ and $s_-$ are of particular interest.
If $D$ is an alternating diagram then $G_{s_+}(D)$ and $G_{s_-}(D)$ 
are the plane graphs first constructed by Tait.

\begin{definition}\label{2.1}\ \
\begin{enumerate}
\item[(i)]
Let $D$ be a diagram of a link and $s$ its Kauffman state. 
We form a graph, $G_s(D)$, associated to $D$ and $s$ as follows.
 Vertices of $G_s(D)$ correspond to circles of $D_s$.
Edges of $G_s(D)$ are in bijection with crossings of $D$ and
an edge connects given 
vertices if the corresponding crossing connects circles of $D_s$ 
corresponding to the vertices\footnote{If $S$ is an enhanced Kauffman 
state of $D$ then, in a similar manner, we associate to $D$ and $S$ the 
graph $G_S(D)$ with signed vertices. Furthermore, we can additionally 
equip $G_S(D)$ with a cyclic ordering of edges at every vertex 
following the ordering of crossings at any circle of $D_s$. 
The sign of each edge is the label of the corresponding crossing.
In short, 
we can assume that $G_S(D)$ is a ribbon (or framed) graph. 
We do not use this additional data in this paper but we plan to 
utilitize this in a sequel paper.}.
\item[(ii)] 
In the language of associated graphs we can state the definition
of adequate diagrams as follows:\ 
the diagram $D$ is $+$-adequate (resp. $-$-adequate) if the graph 
$G_{s_+}(D)$ (resp. $G_{s_-}(D)$) has no loops. 
\end{enumerate}
\end{definition}
In the language of associated graphs we can state the definition 
of adequate diagrams as follows: 
\begin{theorem}\label{2.2} \ \\
Consider a link diagram $D$ of $N$ crossings. Then 
\begin{enumerate}
\item[(+)]
If $D$ is $+$-adequate and $G_{s_+}(D)$ has a cycle 
of odd length, then  the Khovanov homology has $\Z_2$ torsion.
More precisely we show that \\ 
$H_{N-2,N+2|s_+|-4}(D)$ has $\Z_2$ torsion.
\item[(-)]
If $D$ is $-$-adequate and $G_{s_-}(D)$ has a cycle
of odd length, then \\ 
$H_{-N,-N-2|s_-|+4}(D)$ has $\Z_2$ torsion.
\end{enumerate}
\end{theorem}
\begin{proof} 
$(+)$ It suffices to show that the group ${\cal C}_{N-2,N+2|s_+|-4}(D)/
d({\cal C}_{N,N+2|s_+|-4}(D))$ has $2$-torsion.

Consider first the diagram $D$ of the left handed torus knot $T_{2,n}$ 
(Fig.2.1 illustrates the case of $n=5$). The associated graph 
$G_n=G_{s_+}(T_{2,n})$ is an $n$-gon.
\\
\ \\
\centerline{\psfig{figure=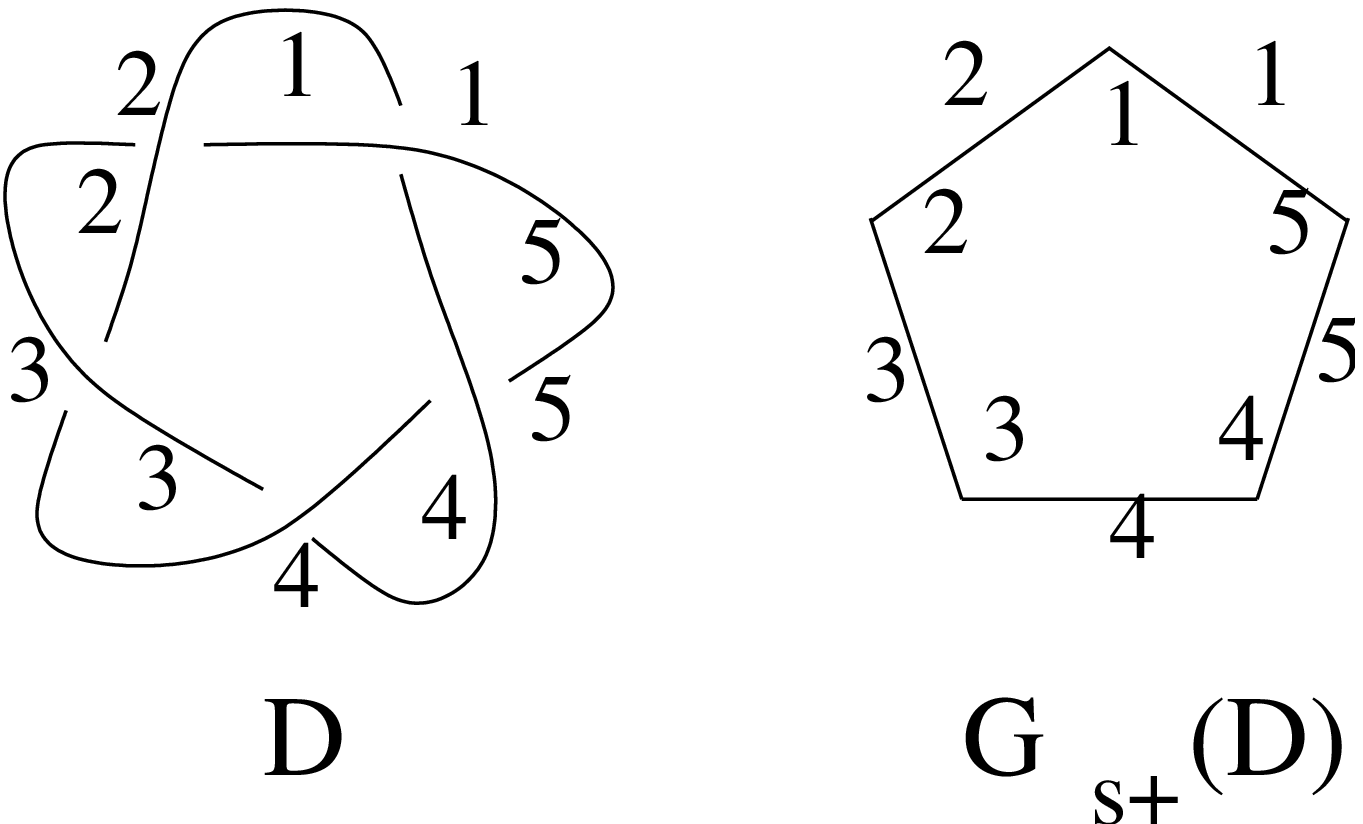,height=4.8cm}}
\begin{center}
Fig. 2.1
\end{center}

 For this diagram we have 
${\cal C}_{n,n+2|s_+|}(D)= \Z$, ${\cal C}_{n,n+2|s_+|-4}(D)=\Z^n$ and 
${\cal C}_{n-2,n+2|s_+|-4}(D) = \Z^n$, where enhanced states generating 
${\cal C}_{n,n+2|s_+|-4}(D)$ have all markers positive and 
exactly one circle (of $D_S$) negative\footnote{In this case $s_+=n$ but we 
keep the general notation so the generalization which follows is natural.}. 
Enhanced states generating
 ${\cal C}_{n-2,n+2|s_+|-4}(D)$ have exactly one negative marker and 
all positive circles of $D_S$. The differential 
$d: {\cal C}_{n,n+2|s_+|-4}(D) \to  {\cal C}_{n-2,n+2|s_+|-4}(D)$ can be 
described by an $n \times n$ circulant matrix (for the ordering of states 
corresponding to the ordering of crossings and regions as in Fig. 2.1)). 

$$\left( \begin{array}{cccccc} 
1 & 1 &  0 & \ldots & 0 & 0\\
0 & 1 & 1 &\ldots & 0 & 0 \\
\ldots & \ldots & \ldots & \ldots & \ldots & \ldots\\ 
0 & 0& \ldots &0 & 1 & 1 \\
1 & 0 & \ldots & 0 & 0 & 1
\end{array} \right)$$

Clearly the determinant of the matrix is equal to 2 (because $n$ is 
odd; for $n$ even the determinant is equal to $0$ because the alternating 
sum of columns gives the zero column). To see this one can 
consider for example the first row expansion\footnote{Because 
the matrix is a circulant one we know furthermore that its eigenvalues are 
equal to $1 + \omega$, where $\omega$ is any $n$'th root of unity 
($\omega^n = 1$), and that $\prod_{\omega^n = 1}(1+\omega)= 0$ for $n$ even 
and $2$ for $n$ odd.}. Therefore the group described by the matrix
is equal to $\Z_2$ (for an even $n$ one would get $\Z$). One more observation 
(which will be used later). The sum of rows of the matrix is equal to 
the row vector ($2,2,2,...,2,2$) but the row vector ($1,1,1,...,1,1$) 
is not an integral linear combination of rows of the matrix. 
In fact the element 
$(1,1,1,...,1,1)$ is the generator of $\Z_2$ group represented by the 
matrix. This can be easily checked because if $S_1,S_2,....S_n$ are 
states freely generating ${\cal C}_{n-2,n+2|s_+|-4}(D)$ then 
relations given by the image of ${\cal C}_{n,n+2|s_+|-4}(D)$ are
$S_2=-S_1, S_3=-S_2=S_1,..., S_1=-S_n=...=-S_1$ thus $S_1+S_2+...+S_n$ is 
the generator of the quotient group ${\cal C}_{n-2,n+2|s_+|-4}(D)/
d({\cal C}_{n,n+2|s_+|-4}(D))= \Z_2$. In fact we have proved that any 
sum of the odd number of states $S_i$ represents the generator of $\Z_2$. 

Now consider the general case in which $G_{s_+}(D)$ is a graph without 
a loop and with an odd polygon. Again, we build a matrix presenting 
the group 
${\cal C}_{N-2,N+2|s_+|-4}(D)/d({\cal C}_{N,N+2|s_+|-4}(D))$ with the
north-west block corresponding to the odd $n$-gon. This block is exactly 
the matrix described previously. Furthermore, the submatrix of the full matrix
 below this block is the zero matrix, as every column 
has exactly two nonzero entries (both equal 
to $1$). This is the case because each edge of the graph (generator) 
has two endpoints (belongs to exactly two relations). If we add  
all rows of the matrix 
 we get the row of all two's. On the other hand 
the row of one's cannot be created, even in the first block. Thus the row 
of all one's representing the sum of all enhanced states in 
${\cal C}_{N-2,N+2|s_+|-4}(D)$ is  $\Z_2$-torsion element in 
the quotient group (presented by the matrix) 
so also in $H_{N-2,N+2|s_+|-4}(D)$. \\
(-) This part follows from the fact that the mirror image of $D$, the diagram 
$\bar D$, satisfies the assumptions of the part (+) of the theorem. 
Therefore the quotient 
${\cal C}_{N-2,N+2|s_+|-4}(\bar D)/d({\cal C}_{N,N+2|s_+|-4}(\bar D))$ has 
$\Z_2$ torsion. Furthermore, the matrix describing the map 
$d: {\cal C}_{-N+2,-N-2|s_-|+4}(D) \to {\cal C}_{-N,-N-2|s_-|+4}(D)$ is
(up to sign of every row) equal to the transpose 
of the matrix describing the map 
$d: {\cal C}_{N,N+2|s_+|-4}(\bar D) \to {\cal C}_{N-2,N+2|s_+|-4}(\bar D)$.
Therefore the torsion of the group 
${\cal C}_{-N,-N-2|s_-|+4}(D)/d({\cal C}_{-N+2,-N-2|s_-|+4}(D)$ is the same 
as the torsion of the group 
${\cal C}_{N-2,N+2|s_+|-4}(\bar D)/d({\cal C}_{N,N+2|s_+|-4}(\bar D))$ 
and, in conclusion, $H_{-N,-N-2|s_-|+4}(D)$ has  $\Z_2$ 
torsion\footnote{Our reasoning reflects a more general fact observed  
 by Khovanov \cite{Kh-1} (see \cite{APS} for the case of $F \times I$)
 that Khovanov homology satisfies ``duality theorem", namely 
$H^{ij}(D) = H_{-i,-j}({\bar D})$. This combined with the 
Universal Coefficients Theorem saying that 
$H^{ij}(D) = H_{ij}(D)/T_{ij}(D) \oplus T_{i-2,j}(D)$, where 
$T_{ij}(D)$ denote the torsion part of $H_{ij}(D)$ gives:
$T_{-N,-N-2|s_-|+4}(D) = T_{N-2,N+2|s_+|-4}(\bar D)$ (notice that 
$|s_-|$ for $D$ equals to $|s_+|$ for $\bar D$).}.

\end{proof}

\begin{remark}\label{2.3}
Notice that the torsion part of the homology, $T_{N-2,N+ 2|s_+|-4}(D)$,
depends only on the graph $G_{s_+}(D)$.
Furthermore if $G_{s_+}(D)$ has no 2-gons then $H_{N-2,N+ 2|s_+|-4}(D)=
{\cal C}_{N-2,N+ 2|s_+|-4}(D)/d(C_{N,N+ 2|s_+|-4}(D))$ and depends only 
on the graph $G_{s_+}(D)$. See a generalization in Remark 3.6
\end{remark}
\section{Diagrams with an even cycle property}\label{3}

If every cycle of the graph $G_{s_+}(D)$ is even (i.e. the graph is 
a bipartite graph) we cannot expect that $H_{N-2,N+|s_+|-4}(D)$  
always has nontrivial torsion. The simplest link diagram 
without an odd cycle in
$G_{s_+}(D)$ is the left handed torus 
link diagram $T_{2,n}$ for $n$ even. As mentioned before, 
in this case ${\cal C}_{n-2,n+2|s_+|-4}(D)/
d({\cal C}_{n,n+2|s_+|-4}(D))= \Z$, and, 
in fact $H_{n-2,n+2|s_+|-4}(D)=\Z$ except $n=2$, i.e. the Hopf link, in 
which case $H_{0,2}(D)=0$.

To find torsion we have to look ``deeper" into the homology. We will 
find a condition for which $H_{N-4,N+2|s_+|-8}(D)$  has  $\Z_2$ torsion, 
where $N$ is the number of crossings of $D$.

Analogously to the odd case, we will start from the left handed torus
link $T_{2,n}$ and associated graph $G_{s_+}(D)$ being an $n$-gon with 
even $n \geq 4$; Fig. 3.1. 
\\
\ \\
\centerline{\psfig{figure=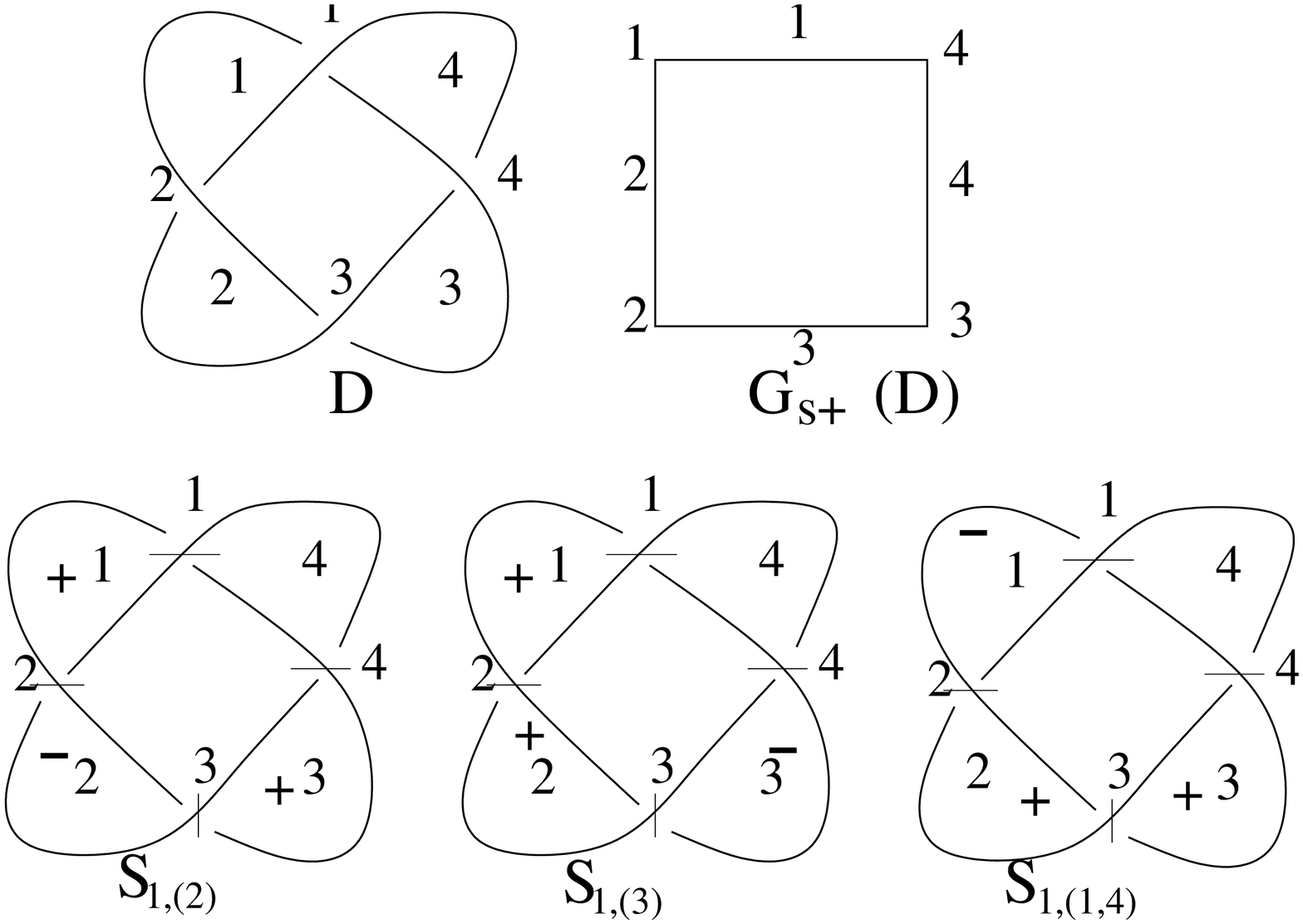,height=5.8cm}}
\begin{center}
Fig. 3.1
\end{center}

\begin{lemma}\label{3.1}\ \\
Let $D$ be the diagram of the left-handed torus link of type 
($2,n$) with $n$ even, $n \geq 4$.\\
Then $H_{n-4,n+2|s_+|-8}(D)= H_{n-4,3n-8}(D) = 
{\cal C}_{n-4,3n-8}(D)/d({\cal C}_{n-2,3n-8}(D)) = \Z_2$.\\ 
Furthermore, every enhanced state from the basis of 
${\cal C}_{n-4,3n-8}(D)$ (or an odd sum of such states) is 
the generator of $\Z_2$.
\end{lemma}
\begin{proof}
We have  $n=|s_+|$. 
The chain group $C_{n-4, 3n-8}(D)= \Z^{\frac{n(n-1)}{2}}$
is freely generated by enhanced states $S_{i,j}$, where exactly 
$i$th and $j$th 
crossings have negative markers, and all the circles of $D_{S_{i,j}}$ are 
positive (crossings of $D$ and circles of $D_{s_+}$ are ordered in Fig. 3.1). 
We have to understand the differential
$d: C_{n-2, 3n-8}(D) \to C_{n-4, 3n-8}(D)$. The chain group 
$C_{n-2, 3n-8}(D) = \Z^{n(n-1)}$ is freely generated by enhanced 
states with $i$th negative marker and one negative circle of $D_S$. 
In our notation we will write $S_{i,(i-1,i)}$ if the negative circle is
obtained by connecting circles $i-1$ and $i$ in $D_{s_+}$ by a negative 
marker. Notation $S_{i,(j)}$ is used if we have $j$th negative circle, 
$j\neq i-1$, $j \neq i$. The states $S_{1,(2)}$, $S_{1,(3)}$ and 
$S_{1,(4,1)}$ are 
shown in Fig. 3.1 ($n=4$ in the figure). The quotient group 
$  {\cal C}_{n-4, 3n-8}(D)/d({\cal C}_{n-2, 3n-8}(D))$ can be presented 
by  a $n(n-1) \times \frac{n(n-1)}{2}$ matrix, $E_n$. One should just 
understand the images of enhanced states of $  {\cal C}_{n-2, 3n-8}(D)$.
In fact, for a fixed crossing $i$ the corresponding $n-1 \times n-1$ 
block is (up to sign of columns\footnote{In the $(n-1)\times (n-1)$ 
block corresponding to the $i$th crossing (i.e. we consider only states 
in which $i$th crossing has a negative marker), the column under the 
generator $S_{i,j}$ of $ {\cal C}_{n-4, 3n-8}$ has $+1$ entries 
if $i<j$ and $-1$ entries if $i>j$.}) 
the circulant matrix discussed in Section 2. 
Our goal is to understand  the matrix $E_n$, to show that 
it represents the group $\Z_2$ and to find natural representatives of 
the generator of the group. 
For $n=4$, $d: \Z^{12} \to \Z^6$ and it is given by: 
$d(S_{1,(2)}) = S_{1,2} + S_{1,3}$, $d(S_{1,(3)}) = S_{1,3} + S_{1,4}$, 
$d(S_{1,(1,4)}) = S_{1,2} + S_{1,4}$, $d(S_{2,(1,2)}) = -S_{2,1} + S_{2,3}$,
 $d(S_{2,(3)}) = S_{2,3} + S_{2,4}$, $d(S_{2,(4)}) = -S_{2,1} + S_{2,4}$, 
$d(S_{3,(1)}) = -S_{3,1} - S_{3,2}$, $d(S_{3,(2,3)}) = -S_{3,2} + S_{3,4}$,
$d(S_{3,(4)}) = -S_{3,1} + S_{3,4}$, $d(S_{4,(1)}) = -S_{4,1} - S_{4,2}$, 
$d(S_{4,(2)}) = -S_{4,2} - S_{4,3}$, $d(S_{4,(4,3)}) = -S_{4,1} - S_{4,3}$,

Therefore $d$ can be described by the $12 \times 6$ matrix. States are 
ordered lexicographically, e.g. $S_{i,j}$ ($i<j$) is before 
$S_{i',j'}$ ($i'<j'$) if $i<i'$ or $i=i'$ and $j<j'$.
$$\left( \begin{array}{cccccc}
1 & 1 & 0 & 0 & 0 & 0\\
0 & 1 & 1 & 0 & 0 & 0 \\
1 & 0 & 1 & 0 & 0 & 0 \\
-1& 0 & 0 & 1 & 0 & 0 \\
0 & 0 & 0 & 1 & 1 & 0 \\
-1& 0 & 0 & 0 & 1 & 0 \\
0 & -1& 0 & -1& 0 & 0 \\
0 & 0 & 0 & -1& 0 & 1 \\
0 & -1& 0 & 0 & 0 & 1 \\
0 & 0 & -1& 0 & -1& 0 \\
0 & 0 & 0 & 0 & -1& -1 \\
0 & 0 & -1& 0 & 0 & -1 
\end{array} \right),$$
In our example the rows correspond to 
$S_{1,(2)},S_{1,(3)},S_{1,(1,4)},S_{2,(1,2)},
S_{2,(3)}, S_{2,(4)}, S_{3,(1)},$ \\
$S_{3,(2,3)},S_{3,(4)},S_{4,(1)},S_{4,(2)},$
and $ S_{4,(4,3)}$, the columns correspond to\\  
$S_{1,2},S_{1,3},S_{1,4},S_{2,3},S_{2,4},S_{3,4}$ 
in this order. 
Notice that the sum of columns of the matrix gives the non-zero 
column of all $\pm 2$ or $0$.
Therefore over $\Z_2$ our matrix represents a nontrivial group. On the 
other hand, over $\Q$, the matrix represent the trivial group. Thus 
over $\Z$ the group represented by the matrix has $\Z_2$ torsion. 
More precisely, we can see that the group is $\Z_2$ as follows: 
The row relations can be expressed as: 
$S_{1,2}=-S_{1,3}=S_{1,4}=-S_{1,2}$, $S_{2,1}= S_{2,3}=-S_{2,4}= -S_{2,1}$, 
$S_{3,1} =- S_{3,2}= -S_{3,4}=-S_{3,1}$ and 
$S_{4,1}= -S_{4,2} = S_{4,3} =-S_{4,1}$. 
$S_{i,j}=S_{j,i}$ in our notation.
In particular, it follows from these 
equalities that the group 
given by the matrix is equal to $\Z_2$ and is generated by 
any basic enhanced state $S_{i,j}$ or the sum of odd number of $S_{i,j}$'s.  

Similar reasoning works for any even $n\geq 4$ (not only $n=4$). 

Furthermore, ${\cal C}_{n-6, 3n-8}=0$, therefore $H_{n-4,3n-8}= 
{\cal C}_{n-4, 3n-8}/d({\cal C}_{n-2, 3n-8}) = \Z_2$.
\end{proof}
We are ready now to use Lemma 3.1 in the general case of an even cycle.
\begin{theorem}\label{3.2}\ \\
Let $D$ be a connected diagram of a link of $N$ crossings such that 
the associated graph $G_{s_+}(D)$ has no loops (i.e. $D$ is $+$-adequate) 
and the graph has an even $n$-cycle with a singular edge 
(i.e. not a part of a 2-gon). 
Then $H_{N-4,N+2|s_+|-8}(D)$  has  $\Z_2$ torsion.
\end{theorem}
\begin{proof}
Consider an ordering of crossings of $D$ such that $e_1,e_2,...,e_n$ are
crossings (edges) of the $n$-cycle. The chain group 
${\cal C}_{N-2,N+2|s_+|-8}(D)$ is freely generated by 
 $N(V-1)$ enhanced
states,$ S_{i,(c)}$, where $N$ is the number of crossings of $D$ (edges
of $G_{s_+}(D)$) and $V=|s_+|$ is the number of circles 
of $D_{s_+}$ (vertices
of $G_{s_+}(D)$). $ S_{i,(c)}$ is the enhanced state in which the crossing
$e_i$ has the negative marker and the circle $c$ of $D_{s_i}$ is negative,
where $s_i$ is the state which has all positive markers except at $e_i$.
The chain group
${\cal C}_{N-4,N+2|s_+|-8}(D)$ is freely generated by enhanced states 
which we can partition into two groups.\\ 
(i) States $S_{i,j}$, where crossings 
$e_i$, $e_j$ have negative markers and corresponding edges of 
$G_{s_+}(D)$ do not form part of a multi-edge 
(i.e. $e_i$ and $e_j$ do not have the 
same endpoints). All circles of the state
$S_{i,j}$ are positive.\\ 
(ii) States $S'_{i,j}$ and $S''_{i,j}$, where crossings
$e_i$, $e_j$ have negative markers and corresponding edges of
$G_{s_+}(D)$ are parts of a multi-edge (i.e. $e_i$, $e_j$ have the same 
endpoints). All but one circle of $S'_{i,j}$ and $S''_{i,j}$ are positive 
and we have two choices for a negative circle leading to 
$S'_{i,j}$ and $S''_{i,j}$, i.e. the crossings $e_i$, $e_j$ touch two circles, 
and we give negative sign to one of them.\\
In our proof we will make the essential use of the assumption that the 
edge (crossing) $e_1$ is a singular edge.

We analyze the matrix presenting the group\\
 ${\cal C}_{N-4,N+2|s_+|-8}(D)/d({\cal C}_{N-2,N+2|s_+|-8}(D))$.

By Lemma 3.1, we understand already the $n(n-1) \times \frac{1}{2}n(n-1)$
block corresponding to the even $n$-cycle.  In this block every column has
$4$ non-zero entries (two $+1$ and two $-1$),
therefore columns of the full matrix corresponding
to states $S_{i,j}$, where $e_i$ and $e_j$ are in the $n$-gon, have zeros
outside our block. We use this property later.

We now analyze another block
represented by rows and columns associated to states having the first
crossing $e_1$ with the negative marker.
This  $(V-1)\times (N-1)$ block has entries equal to $0$ or $1$.
If we add rows in this block we obtain the vector row of two's ($2,2,...,2$),
following from the fact that every edge of $G_{s_+}(D)$  and of $G_{s_1}(D)$
has $2$ endpoints (we use the fact that $D$ is $+$ adequate and $e_1$ is a 
singular edge). Consider
now the bigger submatrix of the full matrix composed of the same rows as
our block but without restriction on columns.
All additional columns are $0$ columns as our row relations
involve only states with negative marker at $e_1$. Thus the sum of these rows
is equal to the row vector ($2,2,...,2,0,...,0$). We will argue now that
the half of this vector, ($1,1,...,1,0,...,0$), is not an integral linear
combination of rows of the full matrix and so represents  $\Z_2$-torsion
element of the group
${\cal C}_{N-4,N+2|s_+|-8}(D)/d({\cal C}_{N-2,N+2|s_+|-8}(D))$.
For simplicity assume that $n=4$ (but the argument holds for any
even $n \geq 4$). Consider the columns indexed by $S_{1,2},S_{1,3},S_{1,4},
S_{2,3},S_{2,4}$ and $S_{3,4}$. The integral linear combination of rows
restricted to this columns cannot give a row with odd number of one's, as
proven in Lemma 3.1. 
In particular we cannot get the row vector ($1,1,1,0,0,0$).
This excludes the row ($1,1,...,1,0,...,0$), as  an integral linear
combination of rows of the full matrix. Therefore the sum of enhanced states
with the marker of $e_1$ negative is  $2$-torsion element in
${\cal C}_{N-4,N+2|s_+|-8}(D)/d({\cal C}_{N-2,N+2|s_+|-8}(D))$ and therefore
in $H_{N-4,N+2|s_+|-8}(D)$.
\end{proof}
Similarly, using duality, we can deal with $-$-adequate diagrams.
\begin{corollary}\label{3.3}\ \\
Let $D$ be a connected, $-$-adequate diagram of a link and the
graph $G_{s_-}(D)$ has an even $n$-cycle, $n \geq 4$, with a singular edge.
Then $H_{-N+2,-N-2|s_-|+8}(D)$  has  $\Z_2$ torsion.
\end{corollary}

\begin{remark}\label{3.4}\ \\
The restriction on $D$ to be a connected diagram is not essential 
(it just simplifies the proof) as for a non-connected diagram, 
$D=D_1\sqcup D_2$ we have ``K\"unneth formula" 
$H_*(D)= H_*(D_1) \otimes H_*(D_2)$ so if any 
of $H_*(D_i)$ has  torsion then $H_*(D)$ has torsion as well.
\end{remark}

We say that a link diagram is doubly $+$-adequate if its graph 
$G_{s_+}(D)$ has no loops and 2-gons. In other words, if a state 
$s$ differs from the state $s_+$ by two markers then $|s| = |s_+| -2 $.
We say that a link diagram is doubly $-$-adequate if its mirror image 
is doubly $+$-adequate.

\begin{corollary}\label{3.5}\ \\
Let $D$ be a connected doubly $+$-adequate diagram of a link of 
$N$ crossings, then either $D$ represents the trivial knot 
or one of the groups
 $H_{N-2,N+2|s_+|-4}(D)$ and $H_{N-4,N+2|s_+|-8}(D)$ has  $\Z_2$ torsion. 
\end{corollary}
\begin{proof}
The associated graph $G_{s_+}(D)$ has no loops and 2-gons.
If $G_{s_+}(D)$ has an odd cycle then 
by Theorem 2.2 $H_{N-2,N+2|s_+|-4}(D)$ has $\Z_2$ torsion. If $G_{s_+}(D)$ 
has an even $n$-cycle, $n \geq 4$ then $H_{N-4,N+2|s_+|-8}(D)$ has 
$\Z_2$ torsion by Theorem 3.2 (every edge of $G_{s_+}(D)$ is a singular edge 
as $G_{s_+}(D)$ has no $2$-gons). Otherwise $G_{s_+}(D)$ is a tree, 
each crossing of $D$ is a nugatory crossing and $D$ 
represents the trivial knot.
\end{proof}

We can generalize and interpret Remark 2.3 as follows.
\begin{remark}\label{3.6} \ \\
 Assume that the associated graph $G_{s_+}(D)$ 
has no $k$-gons, for every $k\leq m$. Then
the torsion part of Khovanov homology, $T_{N-2m,N+ 2|s_+|-4m}(D)$
depends only on the graph $G_{s_+}(D)$.
Furthermore, $H_{N-2m+2,N+ 2|s_+|-4m+4}(D)=
{\cal C}_{N-2m+2,N+ 2|s_+|-4m+4}(D)/d(_{N-2m+4,N+ 2|s_+|-4m+4}(D))$ 
and it depends only on the graph $G_{s_+}(D)$. On a more philosophical 
level
\footnote{In order to be able to recover the full Khovanov homology 
from the graph $G_{s_+}$ we would have to equip the graph with additional 
data: \ ordering of signed edges adjacent to every vertex. This allows us 
to construct a closed surface and the link diagram $D$ on it so that 
$G_{s_+} =G_{s_+}(D)$. The construction imitates the $2$-cell embedding of 
Edmonds (but every vertex corresponds to a circle and signs of edges
 regulate whether 
an edge is added inside or outside of the circle). 
If the surface we obtain is equal to $S^2$ we get the classical Khovanov 
homology. If we get a higher genus surface we have to use \cite{APS} 
theory. This can be  utilitised also to construct  Khovanov homology of 
virtual links (via Kuperberg minimal genus embedding theory \cite{Ku}). 
For example, if the graph $G_{s_+}$ is a loop with adjacent edge(s) 
ordered $e,-e$ then the diagram is composed of a meridian 
and a longitude on the torus.} 
 our observation is related to the fact that if the edge $e_c$ in $G_{s_+}(D)$  
corresponding to a crossing $c$ in $D$ is not a loop then for 
the crossing $c$ the graphs $G_{s_+}(D_0)$ and 
$G_{s_+}(D_{\infty})$ are the graphs obtained from $G_{s_+}(D)$ by 
deleting ($G_{s_+}(D)-e_c$) and contracting ($G_{s_+}(D)/e_c$), 
respectively, the edge $e_c$ (compare Fig. 3.2).
\end{remark}
\centerline{\psfig{figure=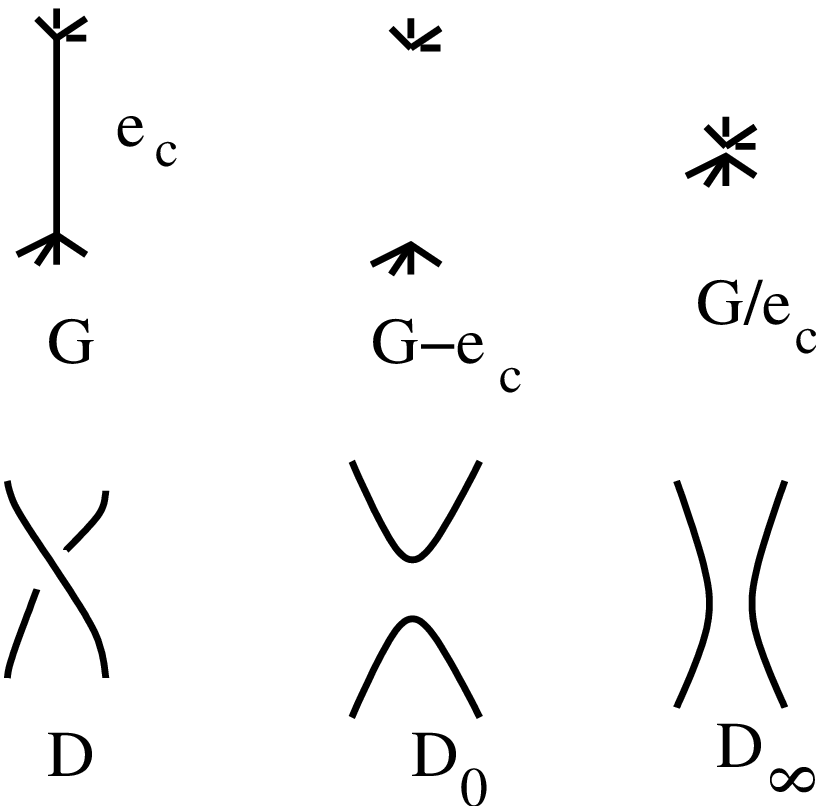,height=4.4cm}}
\begin{center}
Fig. 3.2
\end{center}

\begin{example}\label{3.7}\ \\
 Consider the 2-component alternating link
$6^2_2$  ($\frac{10}{3}$ rational link), with $G_{s_+}(D) = 
G_{s_-}(D)$ being a square with one edge tripled (this is 
a self-dual graph); see Fig. 3.3. 
Corollary 3.5 does not apply to this case but Theorem 3.2
guarantees $\Z_2$ torsion at $H_{2,6}(D)$ and $H_{-4,-6}(D)$. \\
In fact, the KhoHo \cite{Sh-2} computation gives the following 
Khovanov homology\footnote{Tables and programs by Bar-Natan 
and Shumakovitch \cite{BN-3,Sh-2} use 
the version of Khovanov homology for oriented diagrams, 
and the variable 
$q=A^{-2}$, therefore their monomial $q^at^b$ corresponds to the 
free part of the group $H_{i,j}(D;\Z)$ for 
$j=-2b+3w(D)$, $i= -2a +w(D)$ and the monomial $Q^at^b$ corresponds to 
the $\Z_2$ part of the group again with $j=-2b+3w(D)$, $i= -2a +w(D)$. 
KhoHo gives the torsion part of the polynomial for the oriented 
link $6^2_2$, with $w(D)=-6$, as 
$Q^{-6}t^{-1} + Q^{-8}t^{-2} + Q^{-10}t^{-3} + Q^{-12}t^{-4}$.}: 
$H_{6,14}=H_{6,10}=H_{4,10}=\Z$, $H_{2,6}= \Z \oplus \Z_2$, $H_{2,2}=\Z$,
$H_{0,2}= \Z \oplus \Z_2$, $H_{0,-2}= \Z$, $H_{-2,-2} = \Z \oplus \Z_2$,
$H_{-2,-6} = \Z$, $H_{-4,-6}=\Z_2$, $H_{-4,-10}=\Z$, $H_{-6,-10}= H_{-6,-14}=\Z$.
\end{example}

\centerline{\psfig{figure=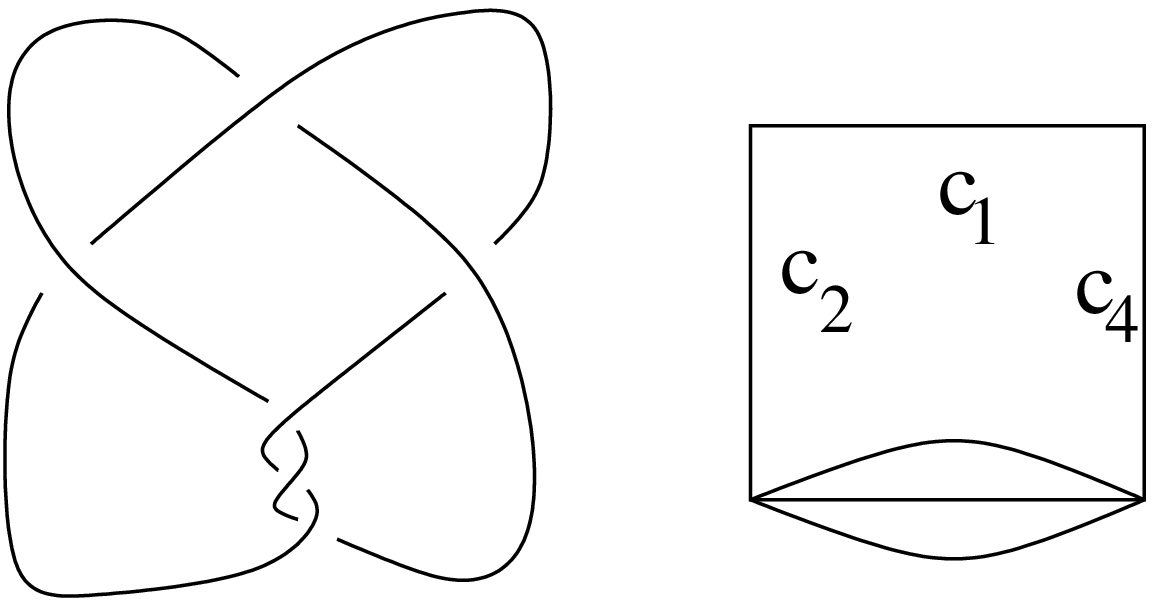,height=4.3cm}}
\begin{center}
Fig. 3.3
\end{center}

\section{Torsion in the Khovanov homology of alternating 
and adequate links}\label{4}

We show in this section how to use technical results from the previous 
sections to prove Shumakovitch's result on torsion in the 
Khovanov homology of alternating links 
and the analogous result for a class of adequate diagrams.

\begin{theorem}[Shumakovitch]\label{4.1} 
The alternating link has torsion free 
Khovanov homology if and only if it is the trivial knot, the Hopf
link or the connected or split sum of copies of them. The nontrivial 
torsion always contains the $\Z_2$ subgroup.
\end{theorem}
The fact that the Khovanov homology of the connected sum of Hopf links 
is a free group, is 
discussed in Section 6 (Corollary 6.6).   

We start with the ``only if" part of the proof by showing 
the following geometric fact.
\begin{lemma}\label{4.2} Assume that $D$ is a link diagram which 
contains a clasp: either $T_{[-2]} = 
\parbox{0.9cm}{\psfig{figure=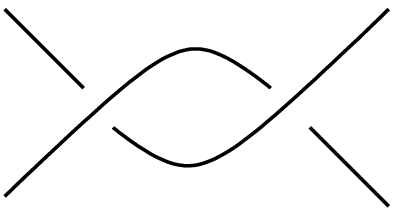,height=0.8cm}}$ \ \ \ \ or 
$T_{[2]} = \parbox{0.9cm}{\psfig{figure=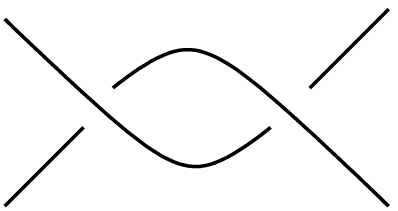,height=0.8cm}}$ 
 \ \ \ \ . 
Assume additionally that the clasp is not a part of the Hopf link summand 
of $D$. Then if the clasp is of 
$T_{[-2]}$ type then the associated graph $G_{s_+}(D)$ has a singular 
edge. If the clasp is of $T_{[2]}$ type then 
the associated graph $G_{s_-}(D)$ has a singular edge. Furthermore the 
singular edge is not a loop.
\end{lemma}
\begin{proof}
Consider the case of the clasp  $T_{[-2]}$, the case of $T_{[2]}$ 
being similar. 
The region bounded by the clasp 
corresponds to the vertex of degree $2$ in $G_{s_+}(D)$. The two edges 
adjacent to this vertex are not loops and they are not singular edges 
only if they share the second 
endpoint as well. In that case our diagram looks like on the 
Fig. 4.1 so it clearly has a Hopf link summand (possibly it is just a 
Hopf link) as the north part is separated by a clasp from the south part 
of the diagram.
\end{proof}
\centerline{\psfig{figure=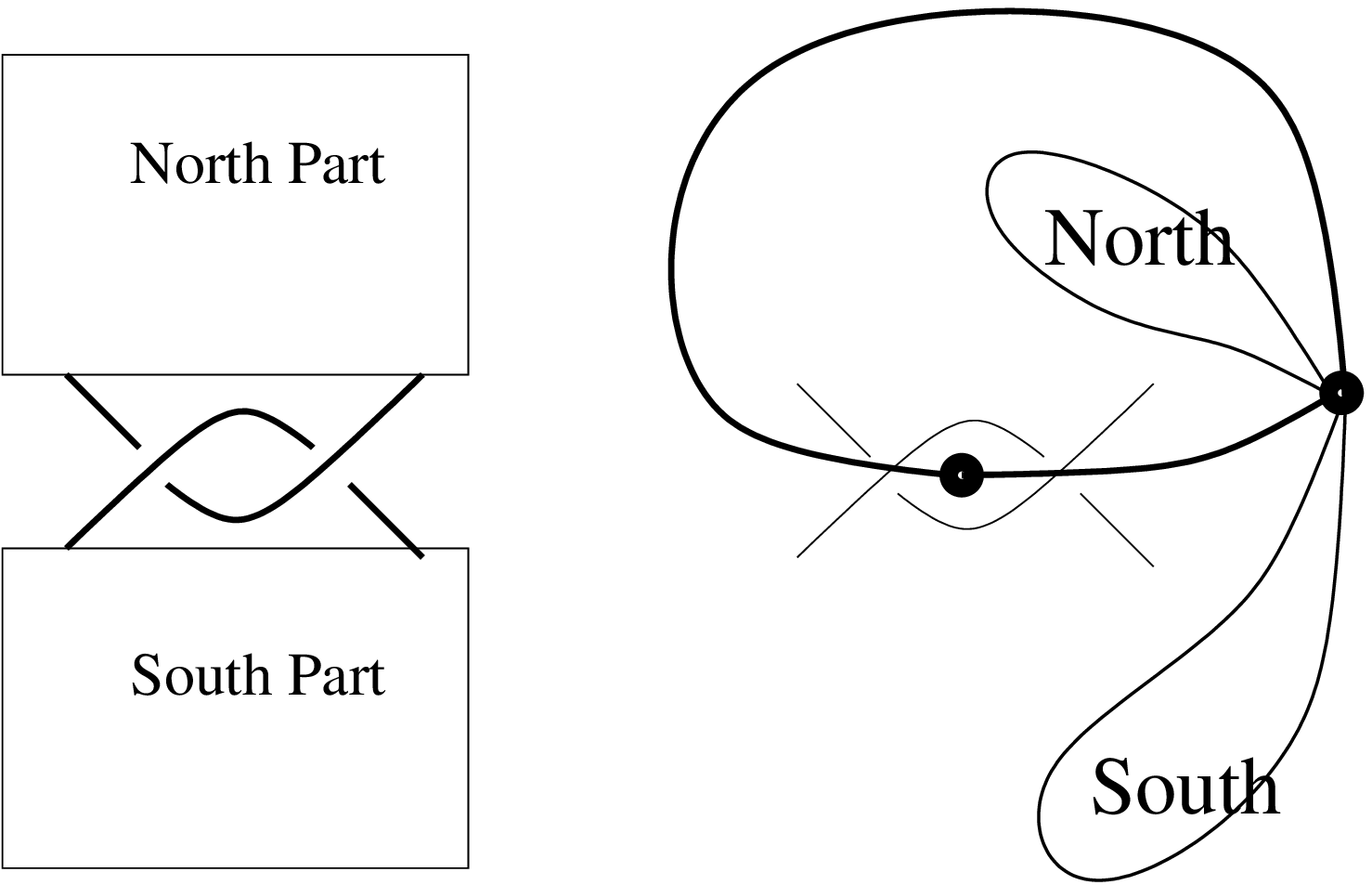,height=3.2cm}}
\begin{center}
Fig. 4.1
\end{center}
\begin{corollary}\label{4.3}
If $D$ is a $+$-adequate diagram (resp. $-$-adequate diagram) with a 
clasp of type $T_{[-2]}$ (resp. $T_{[2]}$), then Khovanov homology 
contains  $\Z_2$-torsion or $T_{[-2]}$ (resp. $T_{[2]}$)) is a part of a 
Hopf link summand of $D$.
\end{corollary}
\begin{proof}
Assume that $T_{[-2]}$ is not a part of Hopf link summand of $D$. By Lemma 
4.2 the graph $G_{s_+}(D)$ has a singular edge. Furthermore,
the graph $G_{s_+}(D)$ has no loops as $D$ is $+$-adequate. If the graph 
has an odd cycle then $H_{N-2,N+2|s_+|-4}(D)$ has $\Z_2$ torsion by Theorem 2.2.
If $G_{s_+}(D)$ is bipartite (i.e. it has only even cycles), then consider 
the cycle containing the singular edge. It is an even cycle of length 
at least 4, so by Theorem 3.2 $H_{N-4,N+2|s_+|-8}(D)$ has $\Z_2$ torsion. 
A similar proof works in $-$-adequate case.
\end{proof}
With this preliminary result we can complete our proof of Theorem 4.1.
 
\begin{proof} First we prove the theorem for non-split, prime 
alternating links. Let $D$ be a diagram of such a link without a 
nugatory crossing. $D$ is an adequate diagram (i.e. it is $+$ and $-$ 
adequate diagram), so it is enough to show that if $G_{s_+}(D)$ (or 
$G_{s_-}(D)$) has a double edge then $D$ can be modified by Tait flypes 
into a diagram with $T_{[-2]}$ (resp. $T_{[2]}$) clasp. This is a standard fact,
justification of which is illustrated in Fig. 4.2\footnote{For alternating 
diagrams, $G_{s_+}(D)$ and $G_{s_-}(D)$ are Tait graphs of $D$. These graphs 
are plane graphs and the only possibilities when multiple edges are not 
``parallel" is if our graphs are not 3-connected (as $D$ is not a split link,
graphs are connected, and because $D$ is a prime link, the graphs are 
2-connected). Tait flype corresponds to the special case of change of the 
graph in its 2-isomorphic class as illustrated in Fig. 4.2.}. 
\\ \ \\
\centerline{\psfig{figure=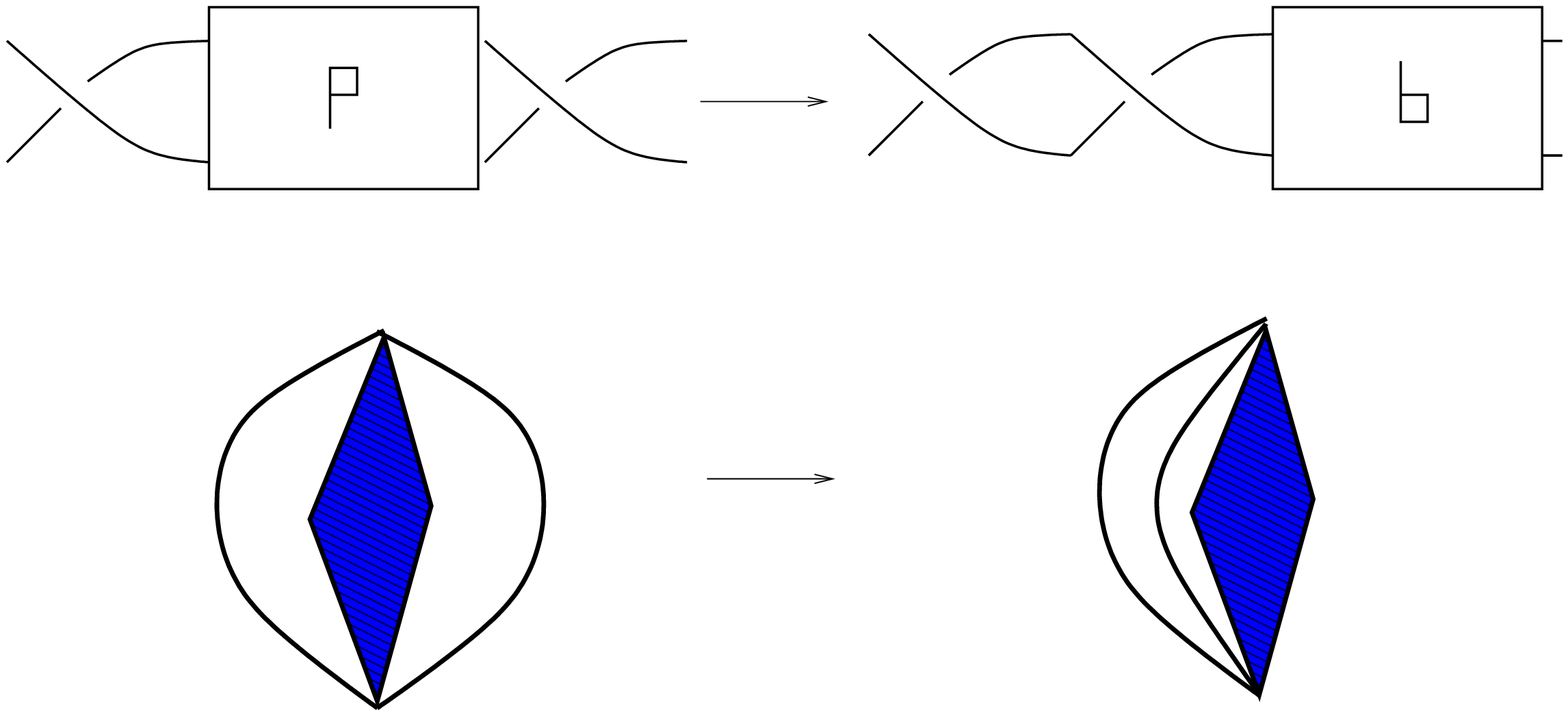,height=4.4cm}}
\begin{center}
Fig. 4.2
\end{center}
If we do not assume that $D$ is a prime link then we use the theorem by 
Menasco \cite{Me} that prime decomposition is visible on the level of 
a diagram. In particular the Tait graphs $G_{s_+}(D)$ and $G_{s_-}(D)$ 
have block structure, where each block (2-connected component) 
corresponds to prime factor of a link. Using the previous results we 
see that the only situation when we didn't find  torsion is if every 
block represents a Hopf link so $D$ represents the sum of Hopf links 
(including the possibility that the graph is just one vertex 
representing the trivial knot).

If we relax condition that $D$ is a non-split link then we use the fact, 
mentioned before, that for $D=D_1\sqcup D_2$, Khovanov homology satisfies 
K\"unneth's formula, $H(D)=H(D_1)\otimes H(D_2)$.
\end{proof}
\begin{example}[The $8_{19}$ knot]\label{4.4} \ \\
The first entry in the knot tables which is not alternating 
is the $(3,4)$ torus knot, $8_{19}$. It is $+$-adequate as it is a positive 
$3$-braid, the closure of $(\sigma_1\sigma_2)^4$. Every positive braid is 
$+$-adequate but its associated graph $G_{s_+}(D)$ 
is composed of 2-gons. Furthermore 
the diagram $D$ of $8_{19}$ is not 
$-$-adequate, Fig. 4.3. 
KhoHo shows that the Khovanov homology of $8_{19}$ has 
 torsion, namely $H_{2,2}=\Z_2$.
This torsion is hidden deeply inside the homology spectrum\footnote{The 
full graded homology group is: $H_{8,14}(D)=H_{8,10}(D)= H_{4,6}(D)=\Z$, 
$H_{2,2}=\Z_2$, $H_{0,2}(D)=H_{2,-2}(D)=H_{0,-2}(D)=H_{-2,-4}(D)= 
H_{-2,-10}(D)=\Z$.}, which 
starts from maximum $H_{8,14}(D)=\Z$ and ends on the minimum 
$H_{-2,-10}(D)=\Z$.\\
\end{example}
\centerline{\psfig{figure=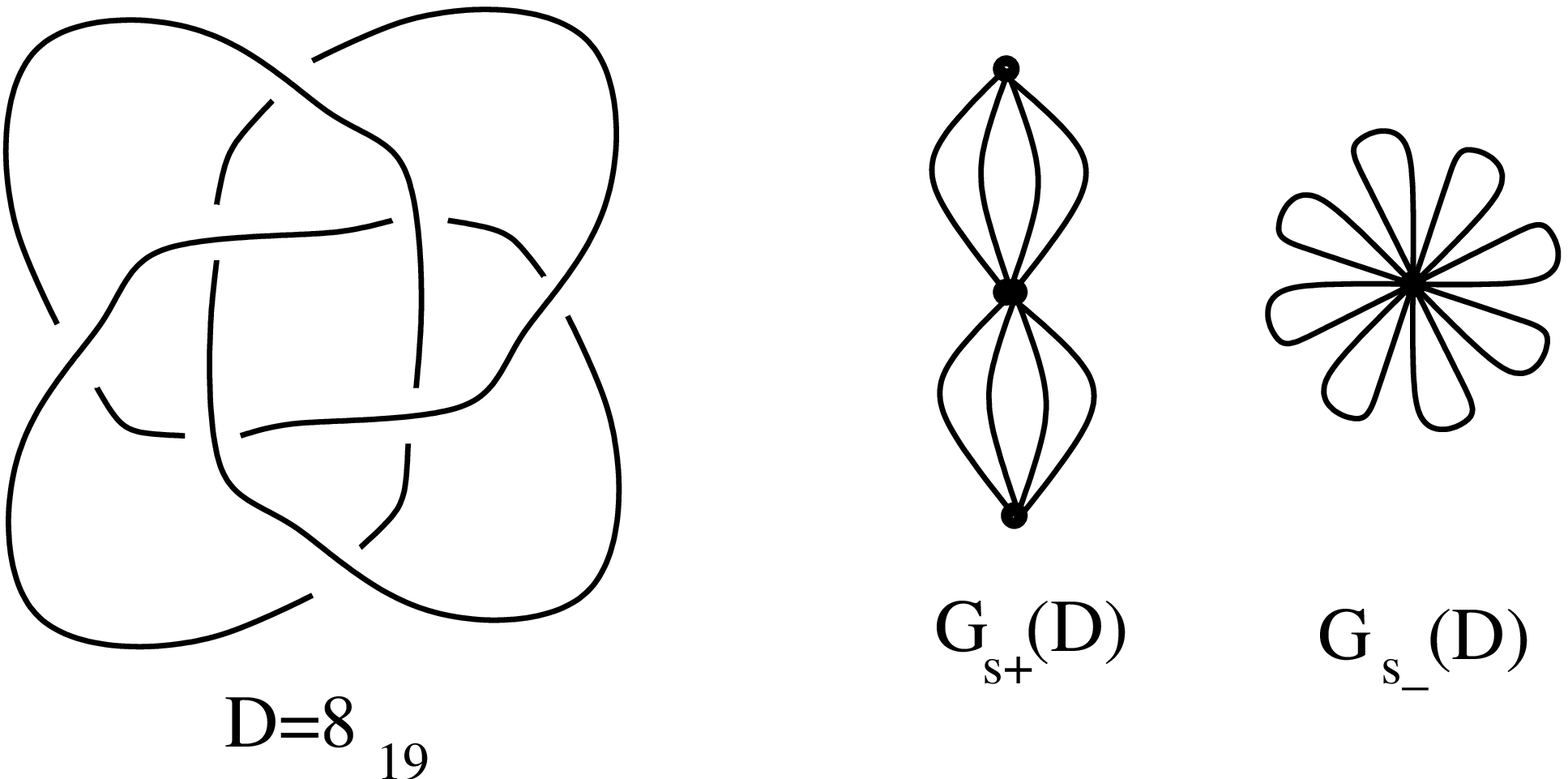,height=5.4cm}}
\begin{center}
Fig. 4.3
\end{center}

\  \\ 

The simplest alternating link which satisfies all conditions of 
Theorem 3.2 except for the existence of a singular edge, is 
the four component alternating link of 8 crossings
$8^4_1$ \cite{Rol}; Fig. 4.4. We know that $H_{**}(8^4_1)$ has torsion 
(by using duality) but Theorem 3.2 does not guarantee torsion in 
$H_{N-4,N+ 2|s_+|-8}(D)= H_{4,8}(D)$, the graph $G_{s_+}(D)$
is a square with every edge doubled; Fig. 4.3. We checked, however 
using KhoHo the torsion part and in fact $T_{4,8}(D)=\Z_2 $. This 
suggests that Theorem 3.2 can be improved\footnote{In \cite{BN-2} the 
figure describes, by mistake, the mirror image of $8^4_1$. 
The full homology is as follows:\ 
$H_{8,16}=\Z$, $H_{8,12}=\Z= H_{6,12}$,
$H_{4,8}=\Z_2 \oplus \Z^4$, $H_{4,4}=\Z$, $H_{2,4}=\Z_2^4$, 
$H_{2,0}=\Z^4$, $ H_{0,0}=\Z^7$, $H_{0,-4}=\Z^6$, 
$H_{-2,-4}=\Z_2^3 \oplus \Z^3$, $H_{-2,-8}=\Z$, $H_{-4,-8}=\Z^3$,
$H_{-4,-12}=\Z^3$, $H_{-6,-12}=\Z_2^3$, 
$H_{-6,-16}=\Z^3$, $H_{-8,-16}=\Z$, $H_{-8,-20}=\Z$. In KhoHo the 
generating polynomials, assuming $w(8^4_1)=-8$, are:
KhPol("8a",21)= $[((q^{18} + q^{16})*t^8 + 3*q^{16}*t^7 + (3*q^{14} + 
3*q^{12})*t^6 + (q^{12} + 3*q^{10})*t^5 + (6*q^{10} + 7*q^8)*t^4 + 
4*q^8*t^3 + (q^6 + 4*q^4)*t^2 + q^2*t + (q^2 + 1))/(q^{20}*t^8)$,\\
$(3*Q^{10}*t^5 + 3*Q^8*t^4 + Q^6*t^3 + 4*Q^2*t + 1)/(Q^{16}*t^6)].$}.
\ \\
\ \\

\centerline{\psfig{figure=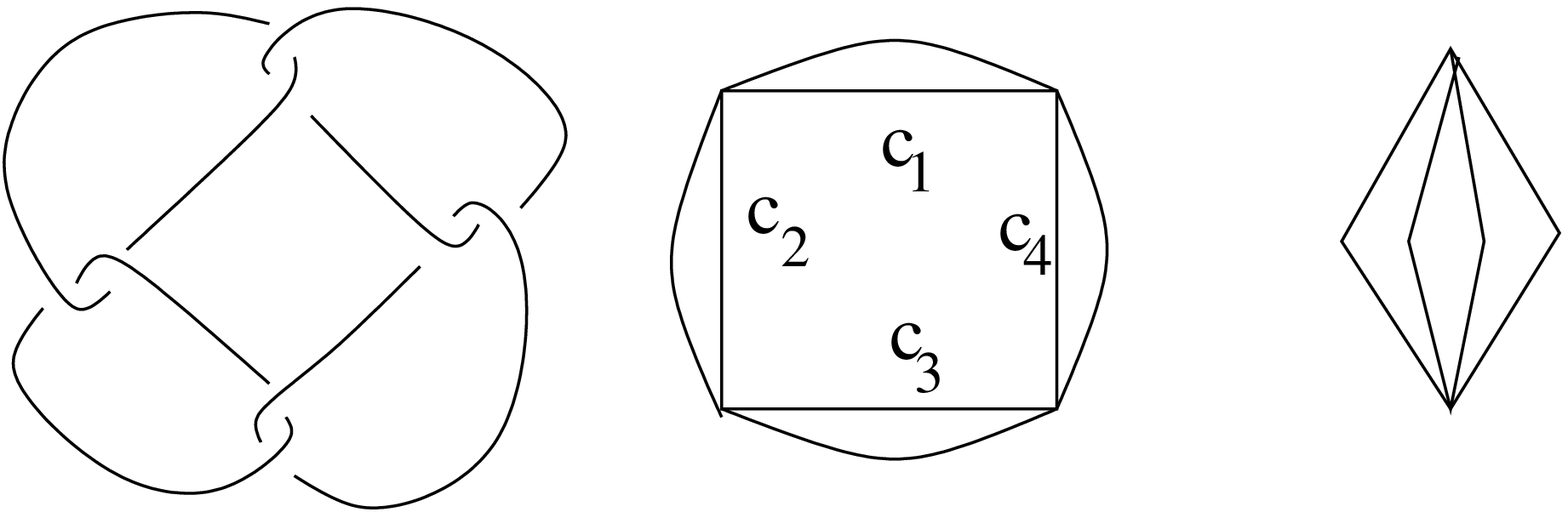,height=4.4cm}}
\begin{center}
Fig. 4.4
\end{center}


\section{Thickness of Khovanov homology 
and almost alternating links}\label{5}

We define, in this section, the notion of an $H$-$k$-thick link diagram and
relate it to $(k-1)$-almost alternating diagrams. In particular we 
give a short proof of Lee's theorem \cite{Lee-1} (conjectured by Khovanov, 
Bar-Natan, and 
Garoufalidis) that alternating non-split links are $H$-$1$-thick 
($H$-thin in Khovanov terminology). 

\begin{definition}\label{5.1}
We say that a link is $k$-almost alternating if it has a diagram 
which becomes alternating after changing $k$ of its crossings.
\end{definition}
As noted in Property 1.4 the ``extreme" terms of Khovanov chain comples 
are $C_{N,N+2|s_+|}(D)= C_{-N,-N-2|s_-|}(D)= \Z$. In the 
following definition 
of a $H$-$(k_1,k_2)$-thick diagram we compare indices of actual 
Khovanov homology of $D$ with lines of slope $2$ going throught the 
points $(N,N+2|s_+|)$ and $(-N,-N-2|s_-|)$.
\begin{definition}\label{5.2}
\begin{enumerate}
\item[(i)] 
We say that a link diagram, $D$ of $N$ crossings is $H$-$(k_1,k_2)$-thick 
if $H_{i,j}(D)=0$ 
with a possible exception of $i$ and $j$ satisfying:
$$ N- 2|s_-| - 4k_2 \leq j-2i \leq 2|s_+| -N + 4k_1.$$

\item[(ii)]
We say that a link diagram of $N$ crossings is 
$H$-$k$-thick\footnote{Possibly, the more appropriate name would be 
$H$-$k$-thin diagram, as the width of Khovanov homology is bounded from 
above by $k$. Khovanov (\cite{Kh-2}, page 7) suggests the term 
homological width; $hw(D)=k$ if homology of $D$ lies on $k$ adjacent 
diagonals (in our terminology, $D$ is $k-1$ thick).}
 if,
it is $H$-$(k_1,k_2)$-thick where $k_1$ and $k_2$ satisfy:
$$k\geq k_1 + k_2 + \frac{1}{2} (|s_+| + |s_-| - N).$$
\item[(iii)] 
We define also $(k_1,k_2)$-thickness (resp. $k$-thickness) 
of Khovanov homology separately for the torsion part 
(we use the notation $TH$-$(k_1,k_2)$-thick 
diagram), and for the free part (we use the notation 
$FH$-$(k_1,k_2)$-thick diagram).
\end{enumerate} 
\end{definition}
Our $FH$-$1$-thick diagram is a $H$-thin diagram 
in \cite{Kh-2,Lee-1,BN-1,Sh-1}.

With the above notation we are able to formulate our main result of 
this section.
\begin{theorem}\label{5.3}\ \\
If the diagram $D_{\infty}= 
D(\parbox{0.9cm}{\psfig{figure=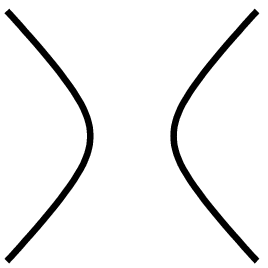,height=0.7cm}})$ 
is $H$-$(k_1(D_{\infty}),
k_2(D_{\infty}))$-thick 
and the diagram 
$D_0 = D(\parbox{0.9cm}{\psfig{figure=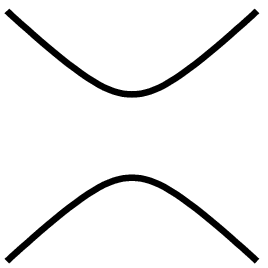,height=0.7cm}})$ 
is $H$-$(k_1(D_0),
k_2(D_0))$-thick, 
then the diagram 
$D_+ = D(\parbox{0.9cm}{\psfig{figure=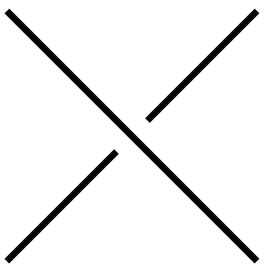,height=0.7cm}})$ 
is $H$-$(k_1(D_+),
k_2(D_+))$-thick 
where \\
$k_1(D_+)= 
\max(k_1(D_{\infty}) + 
\frac{1}{2}(|s_+(D_{\infty})| - 
|s_+(D_+)| +1), 
k_1(D_0))$ and .\\
$k_2(D_+)=
\max(k_2(D_{\infty}),
k_2(D_0)) + \frac{1}{2}(|s_-(D_0)| - |s_-(D_+)| +1))$.\\ 
In particular 
\begin{enumerate}
\item[(i)] if $|s_+(D_+)| - |s_+(D_{\infty})| =1$, as is always 
the case for a $+$-adequate diagram, then 
$k_1(D_+)= \max(k_1(D_{\infty}), k_1(D_0))$,
\item[(ii)] if $|s_-(D_+)| - |s_-(D_0| =1$, as is always
the case for a $-$-adequate diagram, then
$k_2(D_+)= \max(k_2(D_{\infty}), k_2(D_0))$.
\end{enumerate}
\end{theorem}
\begin{proof}
We formulated our definitions so that our proof follows almost 
immediately via the Viro's long exact sequence of Khovanov homology:
$$...\to H_{i+1,j-1}(D_0) \stackrel{\partial}{\to}  
 H_{i+1,j+1}(D_{\infty}) \stackrel{\alpha}{\to} 
H_{i,j}(D_+)
 \stackrel{\beta}{\to}$$
$$H_{i-1,j-1}(D_0) \stackrel{\partial}{\to} H_{i-1,j+1}(D_{\infty})
 \to ...$$
If $0\neq h \in H_{i,j}(D_+)$ then 
either $h=\alpha (h')$ for $0\neq h' \in H_{i+1,j+1}(
D_{\infty})$ or 
 $0\neq \beta(h) \in H_{i-1,j-1}(D_0)$. Thus if
$H_{i,j}(D_+)\neq 0$ 
then either $H_{i+1,j+1}(D_{\infty})\neq 0$ or 
$H_{i-1,j-1}(D_0)\neq 0$.
The first possibility gives the inequalities involving $(j+1)- 2(i+1)$:
$$N(D_{\infty}) - 2|s_-(D_{\infty})| -4k_2(D_{\infty})   
\leq j-2i-1 
\leq 2|s_+(D_{\infty})| - N(D_{\infty}) +4k_1(D_{\infty})$$
which, after observing that $|s_-(D_+)|=|s_-(D_{\infty})|$, leads to:
$$ N(D_+) - 2|s_-(D_+)| - 4k_2(D_{\infty})   \leq j-2i \leq $$\ $$ 
2|s_+(D_+)| - N(D_+) + 4k_1(D_{\infty}) + 2(|s_+(D_{\infty})| - 
|s_+(D_+)| +1).$$
The second possibility gives the inequalities involving $(j-1)- 2(i-1)$:
$$N(D_{0}) - 2|s_-(D_{0})| -4k_2(D_{0})  \leq j-2i+1
\leq 2|s_+(D_{0})| - N(D_{0}) +4k_1(D_{0})$$
which, after observing that $|s_+(D_+)|=|s_+(D_0)|$, leads to:
$$ N(D_+) - 2|s_-(D_+)| - 4k_2(D_{0}) 
-2(|s_-(D_{0})| -|s_-(D_+)| +1) \leq  $$\ $ j-2i \leq \ 
2|s_+(D_+)| - N(D_+) + 4k_1(D_{0}).$

Combining these two cases we obtain the conclusion of Theorem 5.3.
\end{proof}
\begin{corollary}\label{5.4}
If $D$  is an adequate diagram such that, for some crossing of $D$, 
the diagrams $D_0$ and $D_{\infty}$ are $H$-$(k_1,k_2)$-thick (resp. 
$H$-$k$-thick) then $D$ is $H$-$(k_1,k_2)$-thick (resp.  $H$-$k$-thick).
\end{corollary}
\begin{corollary}\label{5.5}
Every alternating non-split diagram without a nugatory 
crossing is $H$-$(0,0)$-thick and $H$-$1$-thick.
\end{corollary}
\begin{proof}
The $H$-$(k_1,k_2)$-thickness in 
Corollary 5.4 follows immediately from Theorem 5.3. To show 
$H$-$k$-thickness we observe additionally that for an adequate 
diagram $D_+$ one has $|s_+(D_0)| + |s_-(D_0)| - N(D_0) =
|s_+(D_+)| + |s_-(D_+)| - N(D_+) = 
|s_+(D_{\infty})| + |s_-(D_{\infty})| - N(D_{\infty})$.

 We prove 
Corollary 5.5 using induction on the number of crossings a slightly
more general statement allowing nugatory crossings.\\ 
(+) If $D$ is an alternating non-split $+$-adequate diagram 
then $H_{i,j}(D)\neq 0$ 
can happen only for $j-2i \leq 2|s_+(D)| - N(D)$.\\
(--) If $D$ is an alternating non-split $-$-adequate diagram 
then $H_{i,j}(D)\neq 0$
can happen only for $N(D) - 2|s_-(D)| \leq j-2i$.\\
If the diagram $D$ from (+) has only nugatory crossings then it represents 
the trivial knot and its nontrivial Khovanov homology are 
$H_{N,3N-2}(D)=H_{N,3N+2}(D)=\Z$. Because $|s_+(D)|=N(D) + 1$ in this case, 
the inequality (+) holds.  In the inductive step we use the property of 
a non-nugatory crossing of a non-split $+$-adequate diagram, namely 
$D_0$ is also an alternating non-split $+$-adequate diagram and inductive 
step follows from Theorem 5.3.\\
Similarly one proves the condition (--). Because the non-split 
alternating diagram without nugatory crossings is an adequate diagram, 
therefore Corollary 5.5 follows from Conditions (+) and (--). 
\end{proof}
 
The conclusion of the theorem is the same if we are interested only in 
the free part of Khovanov homology (or work over a field). In the case 
of the torsion part of the homology we should take into account the 
possibility 
that torsion ``comes" from the free part of the homology, that is 
$H_{i+1,j+1}(D_{\infty})$ may be torsion free but its image under $\alpha$ 
may have torsion element.

\begin{theorem}\label{5.6}
If $T_{i,j}(D_+) \neq 0$ then either \\
(1)\  $T_{i+1,j+1}(D_{\infty}) \neq 0$ or 
$T_{i-1,j-1}(D_{0}) \neq 0$, \\
or\\
(2)\ $FH_{i+1,j+1}(D_{\infty}) \neq 0$ and 
$FH_{i+1,j-1}(D_0) \neq 0$.
\end{theorem}
\begin{proof}
From the long exact sequence of Khovanov homology it follows that the only 
way the torsion is not related to the torsion of $H_{i+1,j+1}(D_{\infty})$ or 
$H_{i-1,j-1}(D_{0})$ is the possibility of torsion created by taking 
the quotient \\ 
$FH_{i+1,j+1}(D_{\infty})/\partial(FH_{i+1,j-1}(D_0))$ and 
in this case both groups $FH_{i+1,j+1}(D_{\infty})$ and 
$FH_{i+1,j-1}(D_0)$ have to be nontrivial.
\end{proof}

\begin{corollary}\label{5.7}
If $D$ is an alternating non-split diagram without a nugatory crossing 
then $D$ is $TH$-$(0,-1)$-thick and $TH$-$0$-thick. In other words if 
$T_{i,j}(D)\neq 0$ then $j-2i = 2|s_+(D)| -N(D)= N(D) -2|s_-(D)|+4$.
\end{corollary}
\begin{proof} We proceed in the same (inductive) manner as in the proof 
of Corollary 5.5, using Theorem 5.7 and Corollary 5.5. In the first step 
of the induction we use the fact that the trivial knot has no torsion 
in Khovanov homology.
\end{proof}

The interest in $H$-thin diagrams was 
motivated by the observation (proved by Lee) that diagrams of non-split 
alternating links are $H$-thin (see Corollary 5.5). Our approach allows 
the straightforward generalization to $k$-almost alternating diagrams.

\begin{corollary}\label{5.8} 
Let $D$ be a non-split $k$-almost alternating diagram without a nugatory 
crossing. Then $D$ is $H$-$(k,k)$-thick and $TH$-$(k,k-1)$-thick.
\end{corollary}
\begin{proof}
The corollary holds for $k=0$ (alternating diagrams) and we use an induction 
on the number of crossings needed to change the diagram $D$ to an alternating 
digram, using Theorem 5.3 in each step.
\end{proof}

We were assuming throughout the section that our diagrams are non-split. 
This assumption was not always necessary. In particular even the split 
alternating diagram without nugatory crossings is $H$-$(0,0)$-thick as 
follows from the following observation.
\begin{lemma}\label{5.9}
If the diagrams $D'$ and $D^{\prime\prime}$ are $H$-$(k_1',k_2')$-thick and 
$H$-$(k_1^{\prime\prime},k_2^{\prime\prime})$-thick, respectively, 
then the diagram $D=D' \sqcup D^{\prime\prime}$ 
is $H$-$(k_1'+k_1^{\prime\prime},k_2'+k_2^{\prime\prime})$-thick.
\end{lemma}
\begin{proof}
Lemma 5.9 follows from the obvious fact that in the split sum 
$D=D' \sqcup D"$ we have 
$N(D) = N(D') + N(D^{\prime\prime})$, $|s_+(D)| = |s_+(D')| + 
|s_+(D^{\prime\prime})|$ and 
$|s_-(D)| = |s_-(D')| + |s_-(D^{\prime\prime})|$.
\end{proof}

Khovanov observed (\cite{Kh-2}, Proposition 7) that adequate 
non-alternating knots are not $H$-$1$-thick. We are able to proof 
the similar result about torsion of adequate non-alternating links.

\begin{theorem}\label{5.10}\ \\
Let $D$ be a connected adequate diagram which  does not represent an 
alternating link and  such that $G_{s_+}(D)$ and 
$G_{s_-}(D)$ have either an odd cycle or an even cycle with a 
singular edge, then $D$ is not $TH$-$0$-thick diagram. More generally,
$D$ is at best $TH$-$\frac{1}{2}(N+2-(|s_+(D)| + |s_-(D)|)$-thick.
\end{theorem}
\begin{proof}
The first part of Theorem 5.10 follows from the second part because by 
Proposition 1.4 (Wu's Lemma), $\frac{1}{2}(N+2-(|s_+(D)| + |s_+(D)|) > 0$ 
for a diagram which is not a connected sum of alternating diagrams. 
By Theorems 2.2, 3.2 and Corollary 3.3, $TH_{i,j}(D)$ is nontrivial 
on slope $2$ diagonals $j-2i = 2|s_+|-N$ and $N- 2|s_-| +4$. The $j$ 
distance between these diagonals is $N- 2|s_-| +4 - (2|s_+|-N) = 
2(N+2-(|s_+(D)| + |s_+(D)|)$, so the theorem follows. 
\end{proof}
\begin{example}\label{5.11}\ \\
Consider the knot $10_{153}$ (in the notation of \cite{Rol}). 
It is an adequate 
non-alternating knot. Its associated graphs $G_{s_+}(10_{153})$ and 
$G_{s_-}(10_{153})$ have triangles (Fig. 5.1) so Theorem 5.10 applies. 
Here $|s_+|=6$, $|s_-|=4$ and by Theorem 2.2, $H_{8,18}(10_{153})$ and
$H_{-10,-14}(10_{153})$ have $\Z_2$ torsion. Thus support of torsion 
requires at least $2$ adjacent diagonals\footnote{Checking \cite{Sh-2},
gives the full torsion of the Khovanov homology of $10_{153}$ as: \ 
$T_{8,18}=T_{4,10}=T_{2,6}=T_{0,6}=T_{-2,-2}=T_{-4,-2}=T_{-6,-6}= 
T_{-10,-14}= \Z_2$.}
\end{example}

\centerline{\psfig{figure=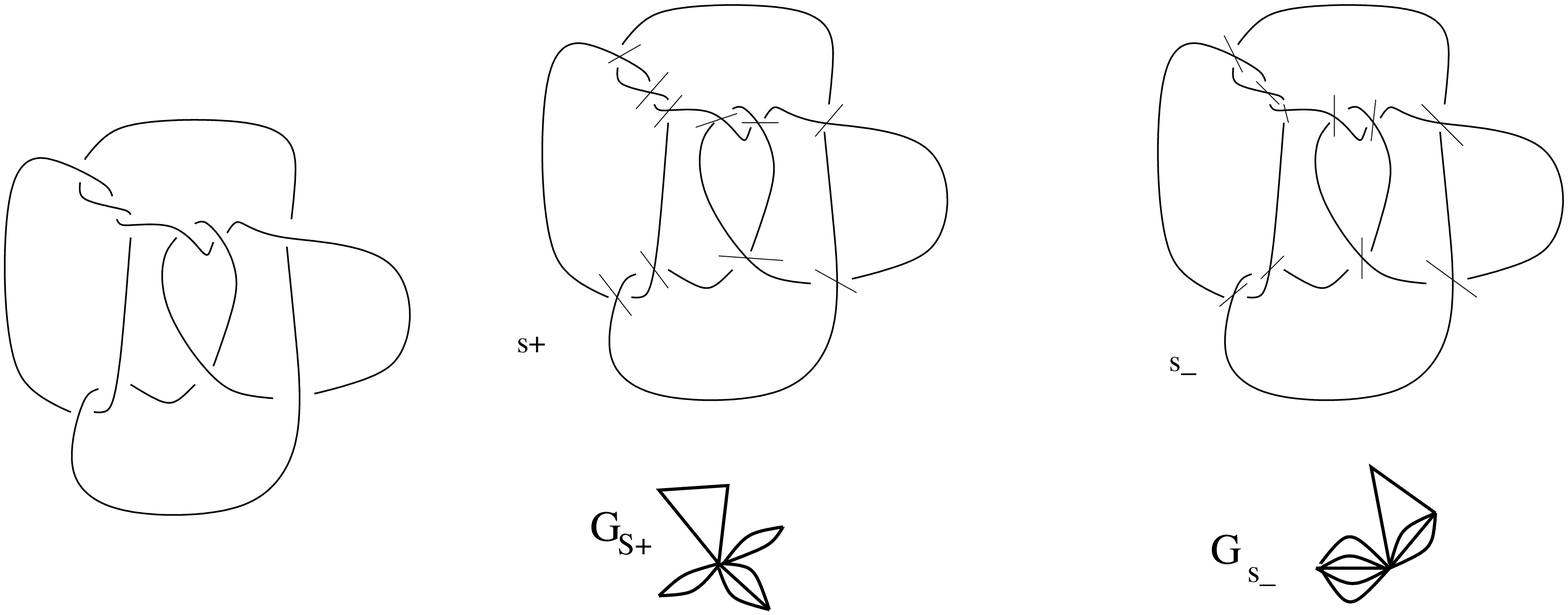,height=5.7cm}}
\begin{center}
Fig. 5.1\end{center}
\begin{corollary}\label{5.12}
Any doubly adequate link which is not an alternating link is not 
$TH$-$0$-thick.
\end{corollary}
 
\section{Hopf link addition}\label{6}
We find, in this section, the structure of the Khovanov homology of 
connected sum of $n$ copies of the Hopf link, as promised in Section 5. 
As a byproduct of our method, we are able to compute Khovanov homology 
of a connected sum of a diagram $D$ and the Hopf link $D_h$, Fig 6.1,  
confirming a conjecture by A.Shumakovitch that the 
Khovanov homology of the connected sum of D with the Hopf link,
is the double of the 
Khovanov homology of $D$.
\\
\ \\
\centerline{\psfig{figure=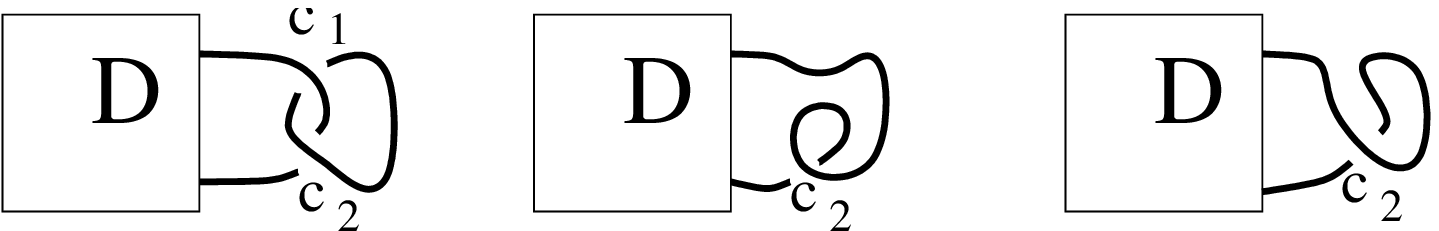,height=2.2cm}}
\centerline{ \ \ \  $D\#D_h$ \ \ \ \ \ \ \ \ \ \ \ \ \ \ \ \ \ \ \ \ 
 $(D\#D_h)_0$ \ \ \ \ \ \ \ \ \ \  \ \ \ \ \ \ \ 
 $(D\#D_h)_{\infty}$  \ \ \ \ } 
\begin{center}
Fig. 6.1
\end{center}

\begin{theorem}\label{6.1}
For every diagram $D$ we have the short exact sequence of Khovanov 
homology\footnote{Theorems 6.1 and 6.2 hold for a diagram $D$ on any surface 
$F$ and for any ring of coefficients $\cal R$ with the restriction 
that for $F=RP^2$ we need $2{\cal R}=0$. In this more general case of 
a manifold being $I$-bundle over a surface, we use definitions and 
setting of \cite{APS}.}
$$0\to H_{i+2,j+4}(D) \stackrel{\alpha_h}{\to} 
H_{i,j}(D\#D_h) \stackrel{\beta_h}{\to} H_{i-2,j-4}(D) \to 0$$
where $\alpha_h$ is given on a state $S$ by Fig.6.2(a)  and $\beta_h$ is a 
projection given by Fig.6.2(b) (and $0$ on other states). The theorem 
holds for any ring of coefficients, $\cal R$, not just ${\cal R}=\Z$.
\end{theorem}
  \ \\
\centerline{\psfig{figure=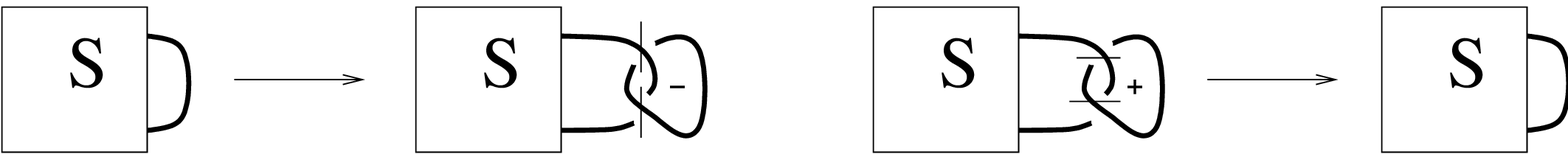,height=1.5cm}}
\ \\
\centerline{ (a) \ \ \ \ \ \ \  \ \ \ \ \ \ \ \ 
$\alpha_h$ \ \ \ \ \ \ \ \ \ \ \ \ \ \ \ \ \ \ \ \ \ \ \ \ \ \ \ \ \ \
 (b) \ \ \ \ \ \ \ \ \ \  \ \ \ \ \ \ \ \ \ \ \ \ \ \ \ \ \ \ \ \
 $\beta_h$  \ \ \ \ \ \ \ \ \ \ \ \ \ } 
\begin{center}
Fig. 6.2
\end{center}
\begin{theorem}\label{6.2}\ \\
The short exact sequence of homology from Theorem 6.1 splits, so we have
$$H_{i,j}(D\#D_h) = H_{i+2,j+4}(D) \oplus H_{i-2,j-4}(D).$$
\end{theorem}
\begin{proof} To prove Theorem 6.1 we 
consider the long exact sequence of the Khovanov homology of the diagram 
$D\#D_h$ with respect to the first crossing of the diagram, $e_1$ (Fig.6.1). 
To simplify the notation we assume that ${\cal R}=\Z$ but our 
proof works for any ring of coefficients.
$$ ... \to H_{i+1,j-1}((D\#D_h)_0) \stackrel{\partial}{\to} 
H_{i+1,j+1}((D\#D_h)_{\infty}) \stackrel{\alpha}{\to} H_{i,j}(D\#D_h) 
\stackrel{\beta}{\to} $$ 
$$H_{i-1,j-1}((D\#D_h)_0) \stackrel{\partial}{\to} 
H_{i-1,j+1}((D\#D_h)_{\infty}) \to ...$$

We show that the homomorphism $\partial$ is the zero map. We 
use the fact that $(D\#D_h)_0$ differs from $D$ by a positive first 
Reidemeister move $R_{+1}$ and that $(D\#D_h)_{\infty}$ 
differs from $D$ by a negative first
Reidemeister move $R_{-1}$; Fig.6.1. 
 We know, see 
\cite{APS} for example, that the chain map \\
$r_{-1}: {\cal C}(D) \to {\cal C}(R_{-1}(D))$ 
given by 
$r_{-1}(\parbox{0.8cm}{\psfig{figure=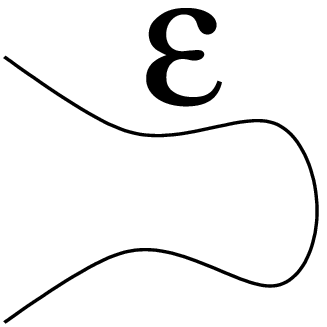,height=0.6cm}}) = 
\parbox{1.2cm}{\psfig{figure=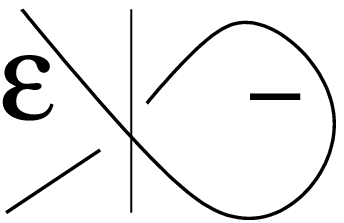,height=0.6cm}}$ 
yields the isomorphism of homology:
$$r_{-1*}: H_{i,j}(D) \to H_{i-1,j-3}(R_{-1}(D))$$ and the chain map 
$\bar{r}_{+1}({\cal C}(R_{+1}(D))={\cal C}((D)$ 
given by the projection with $\bar{r}_{+1}(
\parbox{1.2cm}{\psfig{figure=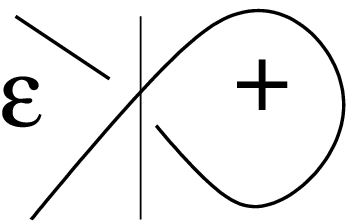,height=0.6cm}})= 
(\ \parbox{1.1cm}{\psfig{figure=kinksmoothe.eps,height=0.6cm}})$ and $0$ 
otherwise,  
induces the isomorphism of homology:
$$\bar r_{+1*}: H_{i+1,j+3}(R_{+1}(D)) \to H_{i,j}(D).$$ 
From these we get immediately that the composition  homomorphism:
$$r_{-1*}^{-1}\partial {\bar r}_{+1*}^{-1}: H_{i,j-4}(D) \to H_{i+2,j+4}(D)$$ 
is the zero map by considering the composition of homomorphisms 
$$ H_{i,j-4}((D)) \stackrel{{\bar r}_{+1*}^{-1}}{\to}
H_{i+1.j-1}((D\#D_h)_{0}) \stackrel{\partial}{\to}
H_{i+1.j+1}((D\#D_h)_{\infty})) \stackrel{r_{-1*}^{-1}}{\to}
H_{i+2,j+4}(D).$$
\end{proof}
Let $h(a,b)(D)$ (resp. $h_{{\cal F}}(a,b)(D)$ for a field ${\cal F}$) 
be the generating polynomial of the free part of $H_{**}(D)$ (resp. 
$H_{**}(D;{\cal F})$),
where $kb^ia^j$  (resp. $k_{{\cal F}}b^ia^j$) 
represents the fact that the free 
part of $H_{i,j}(D)$, $FH_{i,j}(D)= \Z^k$ 
(resp. $H_{i,j}(D;{\cal F})= {\cal F}^k$).

Theorem 6.2 will be proved in several steps.
\begin{lemma}\label{6.3}\ \\
If the module $H_{i-2,j-4}(D;{\cal R})$ is free (e.g. ${\cal R}$ is 
a field) then the sequence from Theorem 6.1 splits and 
 $H_{i,j}((D\#D_h);{\cal R})=H_{i-2,j-4}(D;{\cal R})\oplus 
H_{i+2,j+4}(D;{\cal R})$ or shortly
 $H_{**}(D\#D_h;{\cal R})= 
H_{**}(D;{\cal R})(b^2a^4+ b^{-2}a^{-4})$. 

For the free part we have always $FH_{i,j}((D\#D_h) = 
FH_{i+2,j+4}(D) \oplus FH_{i-2,j-4}(D))$ or in the 
language of generating functions: 
$h(a,b)(D\#D_h)=(b^2a^4+ b^{-2}a^{-4})h(a,b)(D)$.
\end{lemma}
\begin{proof}
The first part of the lemma follows immediately from Theorem 6.1 which holds 
for any ring of coefficients, in particular 
$rank(FH_{i,j}(D\#D_h)=rank(FH_{i+2,j+4}(D)) + rank(FH_{i-2,j-4}(D))$.
\end{proof}
\begin{lemma}\label{6.4}\ \\
There is the exact sequence of $\Z_p$ linear spaces:
$$ 0 \to H_{i+2,j+4}(D)\otimes \Z_p \to H_{i,j}(D\# D_h)\otimes \Z_p
\to H_{i-2,j-4}(D)\otimes \Z_p \to 0 .$$
\end{lemma}
\begin{proof}
Our main tool is the 
universal coefficients theorem
(see, for example, \cite{Ha}; Theorem 3A.3) combined with Lemma 6.3. 
By the second part of Lemma 6.3 it suffices to prove that:
$$T_{i,j}(D\# D_h)\otimes \Z_p = T_{i+2,j+4}(D)\otimes \Z_p \oplus 
T_{i-2,j-4}(D)\otimes \Z_p.$$
From the universal coefficients theorem we have:
$$H_{i,j}(D\#D_h);\Z_p) = 
H_{i,j}(D\#D_h)\otimes \Z_p \oplus  Tor(H_{i-2,j}(D\#D_h), \Z_p)$$ and 
$Tor(H_{i-2,j}(D\#D_h), \Z_p) = T_{i-2,j}(D\#D_h)\otimes \Z_p$ and 
the analogous formulas for the Khovanov homology of $D$. 
Combining this with both parts of 
Lemma 6.3, we obtain:
$$ T_{i,j}(D\# D_h)\otimes \Z_p \oplus T_{i-2,j}(D\# D_h)\otimes \Z_p = $$
$$ 
(T_{i+2,j+4}(D)\otimes \Z_p \oplus T_{i,j+4}(D)\otimes \Z_p) \oplus 
(T_{i-2,j-4}(D)\otimes \Z_p \oplus T_{i-4,j-4}(D)\otimes \Z_p).$$ 
We can express this 
in the language of generating functions assuming that $t(b,a)(D)$ is the 
generating function of dimensions of $T_{i,j}(D)\otimes \Z_p$:
$$(1+b^{-2})t(b,a)(D\# D_h) = (1+b^{-2})(b^2a^4 + b^{-2}a^{-4}t(b,a)(D).$$
Therefore $t(b,a)(D\# D_h) = (b^2a^4 + b^{-2}a^{-4})t(b,a)(D)$ and 
 $dim(T_{i,j}(D\# D_h)\otimes \Z_p) = 
dim(T_{i+2,j+4}(D)\otimes \Z_p + dim(T_{i,j+4}(D)\otimes \Z_p)$.
The lemma follows by observing that 
the short exact sequence with $\Z$ coefficients
leads to the sequence
$$0 \to ker(\alpha_p) \to H_{i+2,j+4}(D)\otimes 
\Z_p \stackrel{\alpha_p}{\to} $$ 
$$H_{i,j}(D\# D_h)\otimes \Z_p
\to H_{i-2,j-4}(D)\otimes \Z_p \to 0 .$$ 
By the previous computation $dim(ker(\alpha_p))=0$ and the 
proof is completed.
\end{proof}
To finish our proof of Theorem 6.2 we only need the following lemma.
\begin{lemma}\label{6.5}
Consider a short exact sequence of finitely generated abelian groups:
$$0\to A \to B \to C \to 0. $$ If for every prime number $p$ we 
have also the exact sequence: 
$$0\to A\otimes \Z_p \to B\otimes \Z_p \to C\otimes \Z_p \to 0 $$ then 
the exact sequence $0\to A \to B \to C \to 0 $ splits and $B= A\oplus C$.
\end{lemma}
\begin{proof} Assume, for contradiction, that the sequence 
$0\to A \stackrel{\alpha}{\to} B \to C \to 0 $ does not split. Then there 
is an element $a \in A$ such that ${\alpha}(a)$ is not $p$-primitive in $B$, 
that is ${\alpha}(a)=pb$ for $b\in B$ and $p$ a prime number and 
$b$ does not lies in the subgroup of $B$ span by ${\alpha}(a)$ 
(to see that such an $a$ exists one 
can use the maximal decomposition of $A$ and $B$ into cyclic subgroups 
(e.g. $A= \Z^{k}\oplus_{p,i} \Z_{p^i}^{k_{p,i}}$)). Now comparing dimensions 
of linear spaces $A\otimes \Z_p, B\otimes \Z_p, C\otimes \Z_p$ 
(e.g. $dim(A\otimes \Z_p)= k+k_{p,1}+k_{p,2}+...$ we see that the sequence 
$0\to A\otimes \Z_p 
\to B\otimes \Z_p \to C\otimes \Z_p \to 0 $ is not exact, a contradiction.
\end{proof}

\begin{corollary}\label{6.6}
For the connected sum of $n$ copies of the Hopf link 
we get\footnote{In the oriented version
(with the linking number equal to $n$, so the writhe number $w=2n$) and
with Bar Natan notation one gets:
$q^{3n}t^n(q+q^{-1})(q^2t + q^{-2}t^{-1})^n$, as computed first by
Shumakovitch.}\\ 
$H_{*,*}(D_h \# ...\# D_h) = h(a,b)(D)=(a^2+a^{-2})(a^4b^2 + a^{-4}b^{-2})^n$
\end{corollary}

\begin{remark}\label{6.7}
Notice that $h(a,b)(D_h)- h(a,b)(OO) = (a^2+a^{-2})(a^4b^2 + a^{-4}b^{-2}) - 
(a^2+a^{-2})^2= b^{-2}a^{-4}(a^2+a^{-2})(1+ba)(1-ba)(1+ba^3)(1-ba^3)$.
This equality may serve as a starting point to formulate a conjecture 
for links, analogous to Bar-Natan-Garoufalidis-Khovanov conjecture 
\cite{Kh-2,Ga},\cite{BN-1} (Conjecture 1), formulated for knots and 
proved for alternating knots by Lee \cite{Lee-1}.
\end{remark}

\section{Reduced Khovanov homology}\label{7}
Most of the results of Sections 5 and 6 can be adjusted to the case of reduced 
Khovanov homology\footnote{Introduced by Khovanov; we follow here 
Shumakovitch's approach adjusted to the framed link version.}. 
We introduce the concept of $H^r$-($k_1,k_2$)-thick diagram and formulate the 
result analogous to Theorem 5.3. The highlight of this section is the 
exact sequence connecting reduced and unreduced Khovanov homology.

Choose a base point, $b$, on a link diagram $D$. Enhanced states, $S(D)$ can 
be decomposed into disjoint union of enhanced states $S_+(D)$ and $S_-(D)$,
where the circle containing the base point is positive, respectively negative.
The Khovanov abelian group 
${\cal C}(D)= {\cal C}_+(D)\oplus {\cal C}_-(D)$ 
where
${\cal C}_+(D)$ is spanned by $S_+(D)$ and ${\cal C}_-(D)$ 
is spanned by $S_-(D)$. ${\cal C}_+(D)$ is a chain subcomplex of 
${\cal C}(D)$. Its homology, $H^{r}(D)$, is called the reduced 
Khovanov homology of $D$, or more precisely of $(D,b)$ (it may 
depends on the component on which the base point lies).
Using the long exact sequence of reduced Khovanov homology we can 
reformulate most of the results of Sections 5 and 6.
\begin{definition}\label{7.1}\ \\
We say that a link diagram, $D$ of $N$ crossings is $H^r$-$(k_1,k_2)$-thick
if $H^r_{i,j}(D)=0$
with a possible exception of $i$ and $j$ satisfying:
$$ N- 2|s_-| - 4k_2 +4 \leq j-2i \leq 2|s_+| -N + 4k_1.$$
\end{definition}

With this definition we have
\begin{theorem}\label{7.2}\ \\
If the diagram $D_{\infty}$ is $H^r$-$(k_1(D_{\infty}),
k_2(D_{\infty}))$-thick
and the diagram
$D_0$ is $H^r$-$(k_1(D_0),
k_2(D_0))$-thick,
then the diagram
$D_+$ is $H^r$-$(k_1(D_+),
k_2(D_+))$-thick
where \\
$k_1(D_+)=
\max(k_1(D_{\infty}) +
\frac{1}{2}(|s_+(D_{\infty})| -
|s_+(D_+)| +1),
k_1(D_0))$ and $k_2(D_+)=
k_2(D_+)=
\max(k_2(D_{\infty}),
k_2(D_0)) + \frac{1}{2}(|s_-(D_0)| - |s_-(D_+)| +1))$.
\\
Every alternating non-split diagram $D$ without a nugatory
crossing is $H^r$-$(0,0)$-thick, and $H^r_{**}(D)$ is torsion free 
\cite{Lee-1,Sh-1}.
\end{theorem}

The graded abelian group ${\cal C}_-(D)= 
\bigoplus_{i,j}{\cal C}_{i,j;-}(D)$ is not a sub-chain complex of ${\cal C}(D)$,
as $d(S)$ is not necessary in ${\cal C}_-(D)$, for $S \in S_-(D)$. 
However the quotient ${\cal C}^-(D)= {\cal C}(D)/{\cal C}_+(D)$ is 
a graded chain complex and as a graded abelian group it can 
be identified with ${\cal C}_-(D)$.
\begin{theorem}\label{7.3}
\begin{enumerate}
\item[(i)] We have the following short exact sequence of chain complexes:
$$ 0 \to {\cal C}_+(D) \stackrel{\phi}{\to}
 {\cal C}(D) \stackrel{\psi}{\to} {\cal C}^-(D)\to 0.$$
\item[(ii)] We have the following long exact sequence of homology:
$$... \to H^r_{i,j}(D) \stackrel{\phi_*}{\to} H_{i,j}(D) 
\stackrel{\psi_*}{\to} H^{\bar r}_{i,j}(D) 
 \stackrel{\partial}{\to} H^r_{i-2,j}(D) \to ...$$ where
$H^{\bar r}_{i,j}(D)$ is the homology of $C^-(D)$. The boundary map 
can be roughly interpreted for a state $S\in S_-(D)$ as $d(S)$ restricted 
to ${\cal C}_+(D)$.
\end{enumerate}
\end{theorem}
 
Applications of Theorem 7.3 will be the topic of a sequel paper, here we
only mention that the group $H^{\bar r}_{**}(D)$ is related to the 
group of the mirror image of $D$, $H^r_{**}(\bar D)$.


 
\section{Acknowledgments}\label{8}

We would like to thank Alexander Shumakovitch for inspiration and 
very helpful discussion.


\noindent
e-mails:\\
\texttt{asaeda@math.uiowa.edu} \\
\texttt{przytyck@gwu.edu}

\end{document}